\documentclass[11pt,leqno]{article}
\usepackage{amsmath,amsthm}
\usepackage{amsfonts}
\usepackage{mathrsfs}
\usepackage{hyperref}
\usepackage{mathtools}
\DeclareRobustCommand{\rchi}{{\mathpalette\irchi\relax}}
\newcommand{\irchi}[1]{\raisebox{\depth}{$#1\chi$}}

\usepackage{amssymb}
\usepackage{epsfig}
\usepackage{subfig}
\usepackage{graphicx}
\usepackage{hyperref}
\usepackage{bm}
\usepackage{booktabs}
\usepackage{color}
\usepackage{float}

\numberwithin{equation}{section} 
\newtheorem{theorem}{\bf Theorem}[section]
\newtheorem{example}{\bf Example}[section]

\newtheorem{remark}{\bf Remark}[section]
\newtheorem{lemma}{\bf Lemma}[section]

\newcommand{\norm}[1]{\left\lVert #1\right\rVert}

\newsavebox{\savepar}

\begin{document}
\title{\bf \LARGE 
Finite Element Analysis for the Chafee-Infante Equation Using Distributed Feedback Control}
\author{
Shishu Pal Singh \thanks{Department of Mathematical Sciences, Rajiv Gandhi Institute of Petroleum Technology. Email: shishups@rgipt.ac.in}  
	\;and
Sudeep Kundu \thanks{Department of Mathematical Sciences, Rajiv Gandhi Institute of Petroleum Technology. Email:sudeep.kundu@rgipt.ac.in }}

\date{\today}
	\maketitle	
	\begin{abstract} In this paper, we propose a \( C^0 \)-conforming finite element method for the Chafee-Infante equation with a finite-parameter feedback control. We establish error analysis for both the state variable and the control variable for the spatially discretized solution. Furthermore, we employ the backward Euler method for time discretization and discuss the stability analysis of the fully discrete scheme. Additionally, we develop error estimates for both the state variable and the control input in the fully discrete setting. Finally, we verify our theoretical conclusions using some numerical experiments.
	\end{abstract}
	

	\section{Inroduction.}
We consider the Chafee-Infante (CI) reaction-diffusion equation with different types of boundary conditions. Specifically, we seek a function $y(x,t)$ for $x \in (0,1)$ and \(t>0\) such that  
 \begin{align}
	\label{1.1}
	\frac{\partial y}{\partial t} - \nu y_{xx} - \gamma y + \delta y^{3} = 0, \\
	y(x,0) = y_{0}(x), \label{1.2}
\end{align}
 with one of the following sets of boundary conditions:
 \begin{subequations}
 	\begin{align}
 		y(0,t) = 0, \quad y_{x}(1,t) = 0, \label{1.3}\\
 		 \quad y(0,t) = 0, \quad y(1,t) = 0, \label{1.14} \\
 		 \quad y_{x}(0,t) = 0, \quad y_{x}(1,t) = 0, \label{1.15}
 	\end{align}
 \end{subequations}
%
where the reaction coefficients $\gamma$ and $\delta$, as well as the diffusion coefficient $\nu$, are positive constants, and $y_{0}(x)$ is a given function.  

Equation \eqref{1.1} is  dissipative reaction-diffusion evolution equation that features a cubic reaction term. Dissipative partial differential equations are widely used in various scientific fields, such as fluid dynamics, chemical reactions, plasma physics, and biological morphogenesis. The complex Ginzburg-Landau (GL) equation is represented by the CI equation, which is the real version. It is also known in the literature as the Newell-Whitehead-Segel equation and the Allen-Cahn equation \cite{Martin2023, inan2020analytical}. The asymptotic behavior of infinite-dimensional dissipative evolution equations can be described using a finite number of degrees of freedom. The CI equation possesses a finite-dimensional global attractor \cite{MR1321480}, which ensures its dissipativity. Additionally, if an inertial manifold exists for a dissipative dynamical system, it can be parameterized using a suitable set of parameters, such as finite volume elements, nodal values, or Fourier modes (see, e.g., \cite{cockburn1997estimating, MR966192, MR1092888, MR972615, JOLLY199038, MR1187736, rosa1997finite, SHVARTSMAN1998637} and references therein).  
The concept of an inertial manifold helps in understanding the long-term behavior of solutions to dissipative dynamical systems. Even though these systems are infinite-dimensional, their inertial manifolds are finite-dimensional. In \cite{MR2016909}, the feedback control law for the CI equation is derived using the concept of the existence of an inertial manifold, which enables the reduction of the system to a finite-dimensional approximation. In contrast, without assuming the existence of an inertial manifold, Azouani and Titi \cite{azouani2013feedback} developed a new type of feedback control strategy that employs various types of determining parameters (e.g., modes, nodes, volume elements, etc.) to control general dissipative dynamical systems.  Furthermore, by applying the feedback control strategy to the system defined by \eqref{1.1}–\eqref{1.2} with zero Neumann boundary conditions \eqref{1.14}, global stabilization in the $L^{2}$- and $H^{1}$-norms has been established.

 For nonlinear parabolic equations and related systems, the finite-parameter feedback control approach has been used for stabilization problems when controllers rely on a finite number of determining parameters. These include the complex Ginzburg–Landau equation, the Navier-Stokes-Voigt equation, the Kuramoto-Sivashinsky (KS) equation, the KdV-Burgers equation, the Burgers equation with nonlocal nonlinearity, the damped nonlinear wave equation, and the Benjamin-Bona-Mahony-Burgers (BBMB) equation among others (see, e.g., \cite{gumus2022finite, KALANTAROVA201740, Varga2018, MR3702656, Ngan2020feedback} and references therein).  Furthermore, the stabilization of the CI equation using boundary feedback control has been demonstrated in \cite{yang2019asymptotic}. Additionally, global exponential stabilization for a one-dimensional nonlinear reaction-diffusion equation has been established in \cite{MR4029808} using a control Lyapunov functional in the $L^{2}$- and $H^{1}$-spatial norms.

In the literature, the numerical study of nonlinear reaction-diffusion equations without control is well-established. Using a nonstandard finite difference method for time discretization and a Galerkin method for spatial discretization, the author in \cite{chin2022analysis} discusses the convergence analysis in the $L^{2}$- and $H^{1}$-norms for the Allen-Cahn equation and demonstrates that the scheme is stable. A similar analysis has been conducted for the real Ginzburg-Landau equation in \cite{chin2021study}.  A semi-analytical Fourier spectral method is presented in \cite{MR3231953}, along with results related to finite difference methods found in \cite{hou2020numerical, inan2018finite, inan2020analytical, rawat2022convergence, MR2679727} for nonlinear reaction-diffusion equations. Additionally, a Haar scale-3 wavelet collocation technique is proposed in \cite{kumar2024application} for the CI equation. Wavelets decompose any function into a summation of basis functions, where each basis function is obtained through dilation and translation of a wavelet function, exhibiting desirable properties of smoothness and locality.

The finite-parameter feedback control algorithm has also been applied in continuous data assimilation, the data assimilation is used in various field such as weather forecasting \cite{azouani2014continuous}, environmental modeling \cite{MR4447757} etc. Continuous data assimilation has been utilized in various studies, including applications to the Navier-Stokes equation \cite{BISWAS2019295, Gesho2016data}, the $2D$ Bénard convection equations \cite{altaf2017Benard, farhat2018assimilation}, general circulation models \cite{desamsetti2019efficient}, and the Allen-Cahn equation \cite{XIA2023511}.  A finite element method has been applied to continuous data assimilation for the $2D$ and $3D$ Navier-Stokes equations in \cite{garcia2020uniform}. The authors discuss two different approaches: a Galerkin method and grad-div stabilization. Moreover, they present a uniform-in-time error analysis between a finite element approximation of the semi-discrete scheme and the reference solution corresponding to the coarse mesh measurements. In \cite{ibdah2020fully}, Ibdah et al. extended these results by incorporating time discretization into the Galerkin spatial approximation scheme from \cite{garcia2020uniform} and established an error analysis for fully discrete schemes in the $L^{2}$- and $H^{1}$-norms.  Further, a fully discrete schemes for data assimilation in the Navier-Stokes equations have been studied in \cite{GARCIAARCHILLA2022114246, MR4122154, LARIOS20191077}. Using finite-parameter feedback control, computational results for the KS equation, CI equation, and other equations have been demonstrated numerically in \cite{MR3702656}.  Using a $C^{0}$-conforming finite element method, error analysis has been conducted for boundary control problems (see, e.g., \cite{sudeep2018burger, sudeep2019BBMB} and references therein). Additionally, stabilization and error analysis of the Neumann boundary problem for the Allen-Cahn equation have been investigated using finite difference analysis in \cite{kundu2024}.

In this paper, we focus on the mixed boundary case \eqref{1.3}. A similar analysis holds for other boundary conditions. Remarks are provided at the end of each section.
Here, we establish an error analysis for the global stabilization of the problem \eqref{1.1} - \eqref{1.3} using finite-parameter feedback control. We apply a $C^ {0}$-conforming finite element method to the CI equation with distributed feedback control. To the best of our knowledge, there is little to no discussion in the literature regarding error analysis for the CI equation with distributed feedback control. Hence, this study aims to provide a rigorous error analysis for both the state variable and the control input.  In this study, we summarize our main results as follows:
\begin{itemize}
    \item  Based on a $C^{0}$-conforming finite element method for spatial discretization, we develop error estimates for the state variable and control input while keeping the time variable continuous.
    \item For time discretization, we use the backward Euler method to analyze an error analysis for the fully discrete scheme. Furthermore, we demonstrate the existence and uniqueness of the fully discrete scheme.
    \item Several numerical experiments are conducted to illustrate stabilization, the order of convergence, and the role of different parameters in the CI equation.
\end{itemize}
In this paper, we adopt the following notations \cite{MR1625845}:
\begin{itemize}
    \item All strongly measurable functions $ v:[0,T]\rightarrow X $ in the space $L^{p}((0,T);X)$ have the norm:
    \begin{align*}
    \norm{v}_{L^{p}((0,T);X)}=
    \begin{cases}
    \left(\int_{0}^{T}\norm{v}_{X}^{p}dt\right)^{\frac{1}{p}}< \infty, &\text{if} \ 1 \leq p<\infty, \\
    \mathop{\mathrm{ess\, sup}} \limits_{0\leq t\leq T} (\norm{v(t)}_{X})< \infty, &\text{if}\  p=\infty.
    \end{cases}
    \end{align*}
Here, \( X \) is a Banach space equipped with the norm \( \| \cdot \|_{X} \), while \( \| \cdot \| \) denotes the standard norm in \( L^{2} \). For simplicity of use, we represent $ L^{p}((0,T);X)$ as $ L^{p}(X).$
    
    \item The Sobolev space is defined as 
    \[
    H^{m}(0,1)=\{v \mid v\in L^{2}(0,1),\ \frac{\partial^{\alpha} v}{\partial x^{\alpha}}\in L^{2}(0,1), \text{ for } \alpha = 1, 2, \ldots, m \}.
    \]
    Additionally, we define $H_{\{0\}}^{1}(0,1):= \{v\in H^{1}(0,1) \mid v(0)=0\}.$ For $ v \in H^{1}(0,1),$ the corresponding norm is given by
    \begin{align*}
    \norm{v}_{1}^{2}:= \norm{v}^{2} +\norm{v_{x}}^{2}, \quad \text{where} \quad \norm{v}^{2}= \int_{0}^{1}v^{2}(x) dx .
        \end{align*}
    
    \item \textbf{Sobolev Embedding:} The following embeddings hold:
    \[
    H^{1}(0,1)\hookrightarrow L^{\infty}(0,1), \quad \text{and} \quad H^{1}(0,1)\hookrightarrow L^{p}(0,1), \quad p\in(2, \infty).
    \]
  
    \item \textbf{Young's Inequality:} For all $ c,d >0$ and $\epsilon>0$,
    \begin{align}
    cd\leq \epsilon\frac{c^{p}}{p}+\frac{1}{\epsilon^{\frac{q}{p}}}\frac{d^{q}}{q},
    \end{align}
    where $ 1<p,q<\infty $ and $ \frac{1}{p} +\frac{1}{q} =1.$
    
\end{itemize}

	The steady-state solution for \eqref{1.1} - \eqref{1.3} is defined as follows:
\begin{align*}
    -\nu y_{xx}^{\infty}- \gamma y^{\infty}+ \delta (y^{\infty})^{3} &= 0,\\
    y^{\infty}(0) = 0, \quad y_{x}^{\infty}(1) &= 0.
\end{align*}
After linearizing equations \eqref{1.1}-\eqref{1.3} about $ y^{\infty}=0 $, we obtain:
\begin{align*}
    \frac{\partial y}{\partial t}-\nu y_{xx}- \gamma y &= 0, \\
    y(x,0) &= y_{0}(x),\\
    y(0,t) &= 0, \quad y_{x}(1,t) = 0.
\end{align*}
Following the proof in \cite{azouani2013feedback} for the Neumann boundary condition \eqref{1.15}, we choose the initial condition 
\[
y_{0}(x) = B_{n} \sin\left(\frac{(2n+1)}{2} \pi x\right),
\]
where $ B_{n} \in \mathbb{R} $. 

Then, the solution takes the form:
\[
y(x,t) = b_{n}(t) \sin\left(\frac{(2n+1)}{2} \pi x\right),
\]
where
$b_{n}(t) = B_{n} e^{\left(-\nu (2n+1)^{2} \frac{\pi^{2}}{4} + \gamma \right)t},$
with the initial condition $ b_{n}(0) = B_{n} $.
Hence, for $ \gamma > 0 $, the term $ b_{n}(t) $ does not converge to zero as $ t \to \infty $ whenever
\(
(2n+1)^{2} < \frac{4\gamma}{\nu \pi^{2}}.
\)
	
The following system represents a distributed feedback control to stabilize the problem \eqref{1.1}-\eqref{1.3} around the steady-state solution $ y^{\infty} \equiv 0 $:
\begin{align}
    \label{1.4}
    \frac{\partial y}{\partial t} - \nu y_{xx} - \gamma y  + \delta y^{3} &= -\mu I_{h}(y), \qquad (x, t) \in (0, 1) \times (0, \infty),\\
    y(0,t) &= 0, \quad y_{x}(1,t) = 0, \qquad t \in (0, \infty),\\
    y(x,0) &= y_{0}(x), \qquad x \in (0,1),
    \label{1.6}
\end{align}
where the identity operator is approximated by the interpolation operator $I_{h}: H^{1}(0,1) \rightarrow L^{2}(0,1) $ with an error of order $h $. We assume the compatibility condition \( y_0(0) = 0, \ y_{0x}(1)=0 \). The parameter $ h $ describes the mesh size, which will be defined later.
\begin{itemize}
	\item Throughout the paper, $ C = C(\delta, \mu, \gamma, c_{p}, \nu) $ denotes a generic positive constant.
\end{itemize}

\begin{itemize}
    \item \cite{MR3702656} The interpolation operator \( I_{h}: H^{1}(0,1) \rightarrow L^{2}(0,1) \) is a general linear operator used to approximate the inclusion map \( i: H^{1} \hookrightarrow L^{2} \), with an associated approximation error of order \( h \). It satisfies the following estimate:
\begin{align}
    \| \phi - I_{h}(\phi) \| \leq c_{p} h \| \phi \|_{1}, \label{1.8}
\end{align}
for all \( \phi \in H^{1}(0,1) \), where \( c_{p} \) is a constant that does not depend on either \( \phi \) or the mesh parameter \( h \). 
\end{itemize}

We define the interpolation operator for nodal values  as follows:

 \textbf{Nodal Values.} For $ \phi \in H^{1}(0,1),$ we define
\begin{align*}
I_{h}(\phi) = \sum_{n=1}^{N} \phi(x_{n}^{*}) \rchi_{J_{n}}(x), \quad x \in [0,1],
\end{align*}
where $ x_{n}^{*} \in J_{n} = [(n-1)\frac{1}{N}, n\frac{1}{N}] $ for $ n = 1, 2, \dots, N $,  and $ \rchi_{J_{n}}(x) $ is the characteristic function of the interval $ J_{n} $. Here, the feedback controllers are represented by the values $ \phi(x_{n}^{*}) $ for $ n = 1, 2, \dots, N $,
 and other interpolation operators are provided in \cite{azouani2013feedback}. The proof of the \eqref{1.8}, using nodal values, is given in \cite{azouani2013feedback} with \(c_{p}=1.\)
 
The rest  of the paper is organized as follows:  Section $2$ discusses regularity results, which plays a crucial role in the later sections. Section $3$ focuses on the finite element methodology for obtaining the semi-discrete solution. It provides a detailed analysis of error estimates associated with the state variable as well as the control input. Section $4$ presents the fully discrete finite element scheme, establishing the existence and uniqueness of the fully discrete solution. Additionally, we provide an error analysis for the state variable in this scheme. Section $ 5,$ concludes with numerical experiments that validate our theoretical findings.

\section{Stabilization.}
In this section, we examine the results of global stabilization for the continuous problem \eqref{1.1}-\eqref{1.3} with the help of a distributed feedback control law. Furthermore, these regularity results are required for error estimates in the state variable and control input.

Let us introduce the weak form of the issue \eqref{1.4} - \eqref{1.6} is to find $ y \in H_{\{0\}}^{1}(0,1) $ such that
\begin{align}
	\label{1.9}
	\left( \frac{\partial y}{\partial t},\chi \right) + \nu (y_{x},\chi_{x}) -  \gamma (y,\chi) + \delta (y^{3},\chi) = -\mu (I_{h}(y),\chi), \quad \forall \ \chi \in H_{\{0\}}^{1}(0,1),
\end{align}
where $ y(x,0) = y_{0}(x) $ is a given function, and $ I_{h}(y) $ is the interpolant operator.

The following assumptions hold throughout this paper:
\begin{align}
\label{2.2}
	\nu\geq \frac{\mu c_{p}^{2}h^{2}}{2}, \quad \text{and}  \quad \mu\geq 2(\gamma+ \nu),
\end{align}
with the decay rate
\begin{align}
\label{2.1}
0< \alpha\leq \frac{1}{2}(\mu- 2\gamma-2\nu),
\end{align}
where \( c_{p} \) comes from \eqref{1.8}.

\begin{remark}
	The  condition \eqref{2.2} impose admissible bounds on the diffusion coefficient \( \nu \) in terms of the feedback gain \( \mu \), parameter \( \gamma \), and mesh size \( h \):
	\[
	\frac{\mu c_p^2 h^2}{2} \leq \nu \leq \frac{\mu}{2} - \gamma.
	\]
	These inequalities ensure that \( \nu \) lies within a stabilizable range consistent with the chosen control parameters. In particular, for a given \( \nu \), the feedback gain \( \mu \) must also be chosen to satisfy
$	\mu \geq 2(\gamma + \nu),$
	ensuring the effectiveness of the control mechanism. See \cite{azouani2014continuous} and \cite{azouani2013feedback} for related results.
\end{remark}

In the following theorem, we show the stabilization in \(H^{2}\)-norm for the CI equation with finite parameter feedback control.

\begin{theorem}
	\label{th2.1}
	Let \( y_{0} \in H^{3}(0, 1)\cap H_{\{0\}}^{1}(0,1) \). Then, the following holds:
	\begin{align*}
		\norm{y}_{2}^{2} + \norm{y_{t}}^{2} + \norm{y_{xt}}^{2} + \norm{y}_{L^{4}}^{4} 
		+ e^{-2\alpha t} \int_{0}^{t} e^{2\alpha s} 
		\left(\norm{y_{t}(s)}_{2}^{2} + \norm{y_{t}(s)}^{2} + \norm{y(s)}_{L^{4}}^{4} \right) ds
		\leq C(\norm{y_{0}}_{3}) e^{-2\alpha t}.
	\end{align*}
\end{theorem}

\begin{proof}
	We refer the reader to Appendix  for the proof.
\end{proof}
\begin{remark}
	Theorem \ref{th2.1} holds in a similar fashion for the Neumann boundary condition \eqref{1.14} and the Dirichlet boundary condition \eqref{1.15}.
\end{remark}

\begin{theorem}
	Let the assumptions \eqref{1.8} and \eqref{2.2} hold, and let the initial data satisfy 
	\[
	y_0 \in H^{3}(0,1) \cap H^{1}_{\{0\}}(0,1).
	\]
	Then the weak solution \( y \) of \eqref{1.9} satisfies the following a priori estimate:
	\begin{align*}
		\norm{y}_{2}^{2} + \norm{y_{t}}^{2} + \norm{y_{xt}}^{2} + \norm{y}_{L^{4}}^{4} 
		+ \int_{0}^{t} \left( \norm{y_{t}(s)}_{2}^{2} + \norm{y_{t}(s)}^{2} + \norm{y(s)}_{L^{4}}^{4} \right) ds 
		\leq C(\norm{y_{0}}_{3}),
	\end{align*}
	for some constant \( C > 0 \) independent of \( t \).
\end{theorem}

\begin{proof}
The proof uses the Bubnov-Galerkin method under assumptions \eqref{1.8} and \eqref{2.2}, along with the \textit{a priori} estimates of Theorem~\ref{th2.1} (with \( \alpha = 0 \)). Uniform bounds are obtained, and the limit is passed using Lions' compactness arguments to establish the result for the weak solution of \eqref{1.9}.

For related results, see \cite{chin2021study, miranville2013existence, MR953967} for Dirichlet, \cite{Georgia2014} for Neumann, and \cite{azouani2013feedback} for uniqueness under Neumann conditions. The method extends to mixed boundary conditions \eqref{1.3}.
\end{proof}
\section{Finite Element Approximation}
In this section, we investigate a semi-discrete Galerkin approximation of problem \eqref{1.4}-\eqref{1.6} while maintaining the time variable continuous, and derive error estimates for both state variable and control input.

For any positive integer \( N \), let \( P = \{ 0 = x_{0} < x_{1} < \ldots < x_{N} = 1 \} \) be a partition of the interval \( (0,1) \) into sub-intervals \( I_{j} = (x_{j-1}, x_{j}) \), with length \( h_{j} = x_{j} - x_{j-1}, \ j = 1, 2, \ldots, N \), and mesh parameter \( h = \max \limits_{1\leq j\leq N} h_{j} \). A finite-dimensional subspace \( S_{h} \subset H_{\{0\}}^{1} \) is defined in the following manner:
\begin{align*}
	S_{h} = \{ v_{h} \in C^{0}(0,1) : v_{h}|_{I_{j}} \in P_{1}(I_{j}), \ \ 0 \leq j \leq N, \ v_{h}(0) = 0 \},
\end{align*}
where \( P_{1}(I_{j}) \) is the space of polynomials of degree \(\leq 1\) on each \( I_{j} \), for \( j = 1, \dots, N \).

Now, the semi-discrete formulation of problem \eqref{1.4}-\eqref{1.6} is to find \( y_{h} \in S_{h} \) such that
\begin{align}
	\label{3.1}
	\left(\frac{\partial y_{h}}{\partial t}, \chi \right) + \nu (y_{hx}, \chi_{x}) - \gamma (y_{h}, \chi) + \delta (y_{h}^{3}, \chi) = -\mu (I_{h}(y_{h}), \chi), \quad \forall \chi \in S_{h},
\end{align}
with the initial condition \( y_{h}(x,0) = y_{0h}(x) \), where \( y_{0h}(x) \) is an approximation of \( y_{0} \), given by \( y_{0h}(x) = P_{h}y_{0}(x) \), where \( P_{h}:H^{1}(0,1) \to S_{h} \) is a projection operator.

Since the space \( S_{h} \) is finite-dimensional and \eqref{3.1} represents a system of ordinary differential equations, a locally unique solution \( y_{h}(t) \) can be obtained using Picard's theorem. In other words, there exists a time \( t = t_{h}^{*} > 0 \) for which the semi-discrete problem \eqref{3.1} admits a unique solution \( y_{h}(t) \) on the interval \( (0, t_{h}^{*}) \). Provided that \( y_{h}(t) \) remains bounded for all \( t > 0 \), a continuation argument is utilized to extend the solution globally in time.

The semi-discrete solution yields the stabilizing result shown below.

\begin{lemma}
	\label{L3.1} 
	Let \( y_{0h} \in L^{2}(0,1) \). Then, we have
	\begin{align*}
	\norm{y_{h}}^{2}+ e^{-2\alpha t} \int_{0}^{t}e^{2\alpha s} \left(\beta \norm{y_{hx}(s)}^{2}+  2\delta \norm{y_{h}(s)}_{L^{4}}^{4} +\beta_{1}\norm{y_{h}(s)}^{2}\right)ds\leq e^{-2\alpha t}\norm{y_{0h}}^{2},
	\end{align*}
	where \( \alpha \) is given in \eqref{2.1}, \( \beta=2\nu- \mu c_{p}^{2}h^{2} \), and \( \beta_{1}=(\mu- 2\gamma-2\nu-2\alpha) \) are positive constants.
\end{lemma}

\begin{proof}
	The proof follows from Theorem \ref{th2.1}. 
\end{proof}

\subsection{Error Estimates}
This subsection determines the semi-discrete scheme's error estimate for the control input and state variable.

We define an auxiliary projection \( \tilde{y}_{h}(t) \in S_{h} \) of \( y(t) \) through the following relation
\begin{align}
	\label{aux}
	(y_{x} - \tilde{y}_{hx}, \chi_{x}) = 0, \quad \forall \chi \in S_{h}.
\end{align}

Let \( \eta := y - \tilde{y}_{h} \) denotes the error in the auxiliary projection \eqref{aux}. For \( \eta \) and \( \eta_{t} \) in the various norms, the estimates are quite standard.

\begin{lemma}
	\label{L4.1}
	There exists a positive constant \( C \) such that for all \( y \in H^{3}(0, 1) \), the following estimates hold
	\begin{align}	
		\norm{\eta} \leq Ch^{2} \norm{y}_{2}, \quad  
		\norm{\eta_{t}} \leq Ch^{2} \norm{y_{t}}_{2}, \quad  
		\text{and} \quad  \norm{\eta_{x}} \leq Ch \norm{y}_{2}.
	\end{align}
\end{lemma}

\begin{proof}
	The proof of the lemma follows from \cite[Chapter 1]{thomee2007galerkin} and \cite{sudeep2018burger}.
\end{proof}	
 We define the error as  
\begin{align*}
	e := y - y_{h} = (y - \tilde{y}_{h}) - (y_{h} - \tilde{y}_{h}) = :\eta - \theta,
\end{align*}
where \( \eta = (y - \tilde{y}_{h}) \) and \( \theta = (y_{h} - \tilde{y}_{h}) \).  
Choose \( \tilde{y}_{h}(0) = y_{0h} \), so that \( \theta(0) = 0 \).

The estimate of \( \eta \) is already given in Lemma \ref{L4.1}; therefore, our primary goal is to determine \( \theta \).  
After deducting \eqref{3.1} from \eqref{1.9} and applying \eqref{aux}, we get
\begin{align}
	 \label{4.7}
	  (\theta_{t},\chi) + \nu (\theta_{x},\chi_{x}) & - \gamma (\theta,\chi) + \delta (\theta^{3},\chi) + 3\delta(\theta^{2}y_{h},\chi) \nonumber \\
	  &= (\eta_{t},\chi) - \gamma(\eta,\chi) + \delta (\eta^{3},\chi) + 3\delta(\eta^{2}y,\chi) - 3\delta(\eta y^{2},\chi) + 3\delta(\theta y_{h}^{2},\chi) \nonumber \\
	  & \quad +\mu (I_{h}(y) - I_{h}(y_{h}),\chi).
\end{align}

\begin{lemma}
	\label{L4.3} Let \( y_{0} \in H^{3}(0,1)\cap H_{\{0\}}^{1}(0,1) \). Then,  
	there exists a decay \(0<2\alpha\leq \beta_{2}\) with \(\mu\geq 3\nu+2\gamma\) and a positive constant \( C =C(\nu, \gamma, \delta, \mu, c_{p}) \), independent of \( h \), such that  
	\begin{align*}
	\norm{\theta}^{2} +e^{-2\alpha t}  \int_{0}^{t} e^{2\alpha s} (\beta \norm{\theta_{x}(s)}^{2} + \delta\norm{\theta(s)}_{L^{4}}^{4}+(\beta_{2}-2\alpha)\norm{\theta(s)}^{2}) ds  
	\leq C h^{4} e^{-2\alpha t} \norm{y_{0}}_{3}^{2},
	\end{align*} 
	where \(\beta_{2}=\mu-3\nu-2\gamma\) and  \(\beta = 2\nu- \mu c_{p}^{2}h^{2}\) are positive constants.
\end{lemma}

\begin{proof}
	Set \( \chi = \theta \) in \eqref{4.7} with \( y_{h}=\theta + \tilde{y_{h}} \) to obtain  
	\begin{align}
	\label{4.8}
	\nonumber \frac{1}{2} \frac{d}{dt} \norm{\theta}^{2} + \nu \norm{\theta_{x}}^{2} - \gamma \norm{\theta}^{2} + \delta \norm{\theta}_{L^{4}}^{4} &= \left( (\eta_{t},\theta) - \gamma (\eta,\theta) \right) + \delta (\eta^{3},\theta) + 3\delta (\eta^{2}y,\theta) - 3\delta (\eta y^{2},\theta) \\&\nonumber  
	\quad + 3\delta (\theta \tilde{y_{h}}^{2},\theta) + 3\delta (\theta^{2} \tilde{y_{h}},\theta) + \mu (I_{h}(y) - I_{h}(y_{h}),\theta), \\&
	= \sum_{i=1}^{4} I_{i}(\theta).
	\end{align}	
	
	On the right side of \eqref{4.8}, the initial term \( I_{1}(\theta) \) yields  
	\begin{align}
	\label{4.9}
	I_{1}(\theta) = \left( (\eta_{t},\theta) - \gamma (\eta,\theta) \right) 
	\leq \frac{\nu}{6} \norm{\theta}^{2} + C\norm{\eta_{t}}^{2} + C(\gamma, \nu)\norm{\eta}^{2}.
	\end{align}			
	On the right-hand side of \eqref{4.8}, the second term \( I_{2}(\theta) \) can be extracted using the Young's inequality to get  
	\begin{align}
	\label{4.10}
	I_{2}(\theta) = \delta (\eta^{3},\theta) 
	\leq C(\delta) \norm{\eta}_{L^{4}}^{4} + \frac{\delta}{4} \norm{\theta}_{L^{4}}^{4}.
	\end{align}	
	With the help of the Young's inequality, the third term \( I_{3}(\theta) \) on the right-hand side of \eqref{4.8} gives
	\begin{align*}
	\nonumber I_{3}(\theta) &= 3\delta (\theta^{2} \tilde{y_{h}},\theta) + 3\delta (\eta^{2} y,\theta) - 3\delta (\eta y^{2},\theta) + 3\delta (\theta \tilde{y_{h}}^{2},\theta), \\&
	\leq C \norm{\tilde{y_{h}}}_{\infty}^{2} \norm{\theta}^{2} + \frac{\delta}{4} \norm{\theta}_{L^{4}}^{4} + \frac{\nu}{6} \norm{\theta}^{2} + C \norm{y}_{\infty}^{2} \norm{\eta}_{L^{4}}^{4} + C \norm{y}_{\infty}^{4} \norm{\eta}^{2} + 3\delta \norm{\tilde{y_{h}}}_{\infty}^{2} \norm{\theta}^{2}.	
	\end{align*}	
	Finally, the last term \( I_{4}(\theta) \) on the right-hand side of \eqref{4.8} is bounded by  
	\begin{align*}
	I_{4}(\theta) &= \mu (I_{h}(y) - I_{h}(y_{h}),\theta)
	=\mu(I_{h}(\eta)-\eta, \theta)+\mu(\eta, \theta)-\mu(I_{h}(\theta)-\theta, \theta)-\mu\norm{\theta}^{2},
	\\[1mm]&\leq C(\nu, \mu, c_{p})h^{2}\norm{\eta}_{1}^{2}+C(\nu, \mu)\norm{\eta}^{2}+ \frac{\nu}{6}\norm{\theta}^{2}+\frac{\mu c_{p}^{2}h^{2}}{2}\norm{\theta_{x}}^{2}+\frac{\mu c_{p}^{2}h^{2}}{2}\norm{\theta}^{2}-\frac{\mu}{2}\norm{\theta}^{2},
	\end{align*}
	
	Substituting these values into \eqref{4.8}, and using the fact that \( H^{1} \xhookrightarrow{} L^{\infty} \) and \( \norm{\tilde{y_{h}}}_{\infty} \leq C\norm{y}_{1} \), it follows from \eqref{2.2} and Theorem \ref{th2.1} with \( \alpha=0 \) that  
	\begin{align*}
	\nonumber \frac{1}{2} \frac{d}{dt} \norm{\theta}^{2} + \frac{1}{2}\beta \norm{\theta_{x}}^{2} + \frac{\delta}{2} \norm{\theta}_{L^{4}}^{4} 
	&\leq \left( \frac{3\nu}{2} +\gamma - \frac{\mu}{2} \right) \norm{\theta}^{2} + C\norm{y}_{1}^{2}\norm{\theta}^{2} \\&
	\quad + C(\norm{\eta_{t}}^{2} + \norm{\eta}^{2} + \norm{\eta}_{L^{4}}^{4}+h^{2}\norm{\eta}_{1}^{2}),
	\end{align*}
	where \(\beta = 2\nu- \mu c_{p}^{2}h^{2}\) is a positive constant.
	
	 Multiplying the above inequality by \( 2e^{2\alpha t} \) and using the Sobolev embedding, we get a modified inequality. We then integrate the resulting inequality  with respect to time from \( 0 \) to \( t \), and apply  Gronwall’s inequality  to obtain
	\begin{align*}
	e^{2\alpha t} \norm{\theta}^{2} &+  \int_{0}^{t} e^{2\alpha s} (\beta \norm{\theta_{x}(s)}^{2} + \delta\norm{\theta(s)}_{L^{4}}^{4}+(\beta_{2}-2\alpha)\norm{\theta(s)}^{2}) ds 
	\\[1mm]&\leq C e^{C\int_{0}^{t}\norm{y(s)}_{1}^{2} ds} \int_{0}^{t} e^{2\alpha s} (\norm{\eta_{t}(s)}^{2} + \norm{\eta(s)}^{2} + \norm{\eta(s)}_{1}^{4}+h^{2}\norm{\eta(s)}_{1}^{2}) ds,
	\end{align*}
	where \(\beta_{2}=\mu-3\nu-2\gamma\geq 0\) and \(\theta(0)=0.\)
	
	Multiplying by \( e^{-2\alpha t} \) and using Lemma \ref{L4.1} and Theorem \ref{th2.1} with \( \alpha=0 \), we obtain  
	\begin{align*}
	e^{2\alpha t} \norm{\theta}^{2} +  \int_{0}^{t} e^{2\alpha s} (\beta \norm{\theta_{x}(s)}^{2} + \delta\norm{\theta(s)}_{L^{4}}^{4}+(\beta_{2}-2\alpha)\norm{\theta(s)}^{2}) ds  
	\leq e^{-2\alpha t} C h^{4} \norm{y_{0}}_{3}^{2}.
	\end{align*}
	This completes the proof.
\end{proof}



\begin{remark}
	\label{r3.2}
	To find the error corresponding to the Neumann boundary conditions, we modify the auxiliary projection \( \tilde{y_{h}} \in S_{h} \) of \( y(t) \) in the following form:
	\begin{align*}
		(y_{x}-\tilde{y}_{hx},\chi_{x})+ \lambda (y-\tilde{y}_{h},\chi)=0, \quad \chi \in S_{h},
	\end{align*}
	where \( \lambda \) is a positive number. The Lax-Milgram Lemma explicitly implies the existence of a unique \( \tilde{y_{h}}(t) \). For more details on this type of projection, see \cite{doi:10.1137/0710062}.

Then equation \eqref{4.7} takes the following form: 
\begin{align*}
\nonumber (\theta_{t},\chi) + \nu (\theta_{x},\chi_{x}) & - \gamma (\theta,\chi) + \delta (\theta^{3},\chi) + 3\delta (\theta^{2}y_{h},\chi)\\
& = (\eta_{t},\chi) - (\lambda \nu+ \gamma)(\eta,\chi) + \delta (\eta^{3},\chi) + 3\delta (\eta^{2}y,\chi) - 3\delta (\eta y^{2},\chi) + 3\delta (\theta y_{h}^{2},\chi)\\
& \quad +\mu (I_{h}(y) - I_{h}(y_{h}),\chi).
\end{align*}

For the proof of the error analysis, we proceed similarly to the case of mixed boundary conditions. In the case of Dirichlet boundary conditions, the error analysis proof follows in the same manner as in the mixed boundary case.
\end{remark}
\begin{lemma}
	\label{L4.4}
	Let \( y_{0} \in H^{3}(0, 1)\cap H_{\{0\}}^{1}(0,1) \). Then, there is a positive constant \( C \) independent of \( h \) such that
	\begin{align*}
	\nu \norm{\theta_{x}}^{2} + \frac{1}{2} e^{-2\alpha t} \int_{0}^{t} e^{2\alpha s} \norm{\theta_{t}(s)}^{2} ds + \frac{\delta}{2} \norm{\theta}_{L^{4}}^{4} 
	\leq C(\norm{y_{0}}_{3}) e^{-2\alpha t} h^{4} e^{C \norm{y_{0}}_{1}^{2}}.
	\end{align*}
\end{lemma}
\begin{proof}
	Setting \( \chi=\theta_{t} \) in \eqref{4.7}, we obtain
\begin{align}
\label{4.15}
\norm{\theta_{t}}^{2} + \frac{\nu}{2} \frac{d}{dt} \norm{\theta_{x}}^{2} - \gamma (\theta, \theta_{t}) + \frac{\delta}{4} \frac{d}{dt} \norm{\theta}_{L^{4}}^{4} = \sum_{i=1}^{4} I_{i}(\theta_{t}).
\end{align}
On the right-hand side of \eqref{4.15}, the first term \( I_{1}(\theta_{t}) \) yields
\begin{align*}
I_{1}(\theta_{t}) = \left( (\eta_{t} ,\theta_{t}) - \gamma(\eta , \theta_{t}) \right) 
\leq \frac{1}{14} \norm{\theta_{t}}^{2} + C\norm{\eta_{t}}^{2} + C(\gamma) \norm{\eta}^{2}.
\end{align*}
Using Young's inequality, the second term \( I_{2}(\theta_{t}) \) gives
\begin{align*}
I_{2}(\theta_{t}) = \delta(\eta^{3} , \theta_{t}) \leq C(\delta) \norm{\eta}_{L^{6}}^{6} + \frac{1}{14} \norm{\theta_{t}}^{2}.
\end{align*}
On the right side of \eqref{4.15}, the third term \( I_{3}(\theta_{t}) \) is bounded by substituting \( y_{h} = \tilde{y_{h}} + \theta \) 
\begin{align*}
I_{3}(\theta_{t}) &= 3\delta(\eta^{2} y , \theta_{t}) - 3\delta(\eta y^{2} , \theta_{t}) + 3\delta (\theta^{2} \tilde{y_{h}} , \theta_{t} ) + 3\delta (\theta(\tilde{y_{h}})^{2} , \theta_{t}),
\\&\leq C(\delta) \norm{\theta}_{L^{4}}^{4} \norm{\tilde{y_{h}}}_{\infty}^{2} + C(\delta) \norm{\eta}_{L^{4}}^{4} \norm{y}_{\infty}^{2} 
+ C(\delta) \norm{\eta}^{2} \norm{y}_{\infty}^{4} + C(\delta) \norm{\theta}^{2} \norm{\tilde{y_{h}}}_{\infty}^{4} + \frac{2}{7} \norm{\theta_{t}}^{2}.
\end{align*}
Finally, the last term on the right-hand side, \( I_{4}(\theta_{t}) \), can be expressed as
\begin{align*}
I_{4}(\theta_{t}) &= \mu (I_{h}(y)- I_{h}(y_{h}) , \theta_{t})
\leq C(\norm{\eta}^{2}+h^{2}\norm{\eta}_{1}^{2})+Ch^2\norm{\theta}_{1}^{2} + \frac{1}{14} \norm{\theta_{t}}^{2}.
\end{align*}
Hence, from \eqref{4.15}, using Theorem \ref{th2.1} with \( \norm{\tilde{y_{h}}}_{\infty} \leq C \norm{y}_{1} \), we obtain
\begin{align*}
\frac{\nu}{2} \frac{d}{dt} \norm{\theta_{x}}^{2} + \frac{1}{2} \norm{\theta_{t}}^{2} + \frac{\delta}{4} \frac{d}{dt} \norm{\theta}_{L^{4}}^{4}
&\leq \gamma (\theta, \theta_{t}) + C(\norm{\eta_{t}}^{2} + \norm{\eta}^{2} + \norm{\eta}_{L^{4}}^{4} + \norm{\eta}_{L^{6}}^{6}+h^{2}\norm{\eta}_{1}^{2}) \\ 
&\quad + C \norm{\theta}_{L^{4}}^{4} \norm{y}_{1}^{2} + C \norm{\theta}^{2} \norm{y}_{1}^{4} + \frac{7\mu^{2}c_{0}^{2}}{2} \norm{\theta}^{2}+Ch^2\norm{\theta}_{1}^{2}.
\end{align*}
Using the Young's inequality 
\[\gamma (\theta, \theta_{t})\leq \frac{1}{4}\norm{\theta_{t}}^{2}+ \gamma^{2}\norm{\theta}^{2},
\]
therefore, with the help of the Sobolev embedding along with Lemma \ref{L4.1}, we arrive at
\begin{align*}
\frac{\nu}{2} \frac{d}{dt} \norm{\theta_{x}}^{2} + \frac{1}{4} \norm{\theta_{t}}^{2} + \frac{\delta}{4} \frac{d}{dt} \norm{\theta}_{L^{4}}^{4}
&\leq Ch^{4} (\norm{y_{t}}_{2}^{2} + \norm{y}_{2}^{2} + \norm{y}_{2}^{4} + h^{2} \norm{y}_{2}^{6}+h^{2}\norm{y}_{2}^{2}) \\
&\quad + C \norm{\theta}_{L^{4}}^{4} \norm{y}_{1}^{2} + C \norm{\theta}^{2} \norm{y}_{1}^{4} + (\gamma^2 + \frac{7\mu^{2}c_{0}^{2}}{2}) \norm{\theta}^{2}+Ch^2\norm{\theta}_{1}^{2},
\end{align*}
where \(\norm{\eta}_{L^{6}}^{6}\leq C \norm{\eta}_{1}^{6}\leq Ch^{6}\norm{y}_{2}^{6}\).

%
Multiplying the above inequality by \( 2e^{2\alpha t} \) and then integrating from \( 0 \) to \( t \), we obtain
\begin{align*}
\nu \norm{\theta_{x}}^{2} &+ \frac{1}{2} e^{-2\alpha t} \int_{0}^{t} e^{2\alpha s} \norm{\theta_{t}(s)}^{2} ds + \frac{\delta}{2} \norm{\theta}_{L^{4}}^{4} \\
&\leq Ch^{4} e^{-2\alpha t} \int_{0}^{t} e^{2\alpha s} (\norm{y_{t}(s)}_{2}^{2} + \norm{y(s)}_{2}^{2} + \norm{y(s)}_{2}^{4} + h^{2} \norm{y(s)}_{2}^{6}+h^{2}\norm{y(s)}_{2}^{2}) ds \\ 
&\quad + C e^{-2\alpha t} \int_{0}^{t} e^{2\alpha s} ( \norm{\theta(s)}_{L^{4}}^{4} + \norm{\theta(s)}^{2} \norm{y(s)}_{1}^{2}) \norm{y(s)}_{1}^{2} ds \\ 
&\quad + (2\gamma^2 + 7\mu^{2}c^{2}) e^{-2\alpha t} \int_{0}^{t} e^{2\alpha s} \norm{\theta(s)}^{2} ds + C e^{-2\alpha t} \int_{0}^{t} e^{2\alpha s} (h^2\norm{\theta(s)}_{1}^{2}+\norm{\theta(s)}_{L^{4}}^{4}) ds.
\end{align*}

Therefore, using  Lemmas \ref{L4.1} - \ref{L4.3} and Theorem \ref{th2.1} with \( \alpha = 0 \), the proof is completed.
\end{proof}

The proof of the following theorem follows from Lemmas \ref{L4.3} - \ref{L4.4} with the help of the triangle inequality.
\begin{theorem}
	\label{th3.1}
   Suppose the hypotheses of Lemmas \ref{L4.3} - \ref{L4.4} hold. Then, there exists a constant $ C $ independent of $ h $ such that 
   \begin{align*}
   	\norm{y-y_{h}}_{i} \leq C h^{2-i}, \quad i=0,1.
   \end{align*}
\end{theorem}

In the following theorem, we discuss the error associated with the control input.

\begin{theorem}
	\label{th3.2}
	Assume that the hypotheses of Theorem \ref{th3.1} are satisfied. Then, the following estimate holds:
	\begin{align*}
	\norm{I_{h}(y)-I_{h}(y_{h})} \leq Ch^{2},
	\end{align*}
	where $ C $ is a positive constant independent of $ h. $
\end{theorem}

\begin{proof}
	Since the interpolation operator is linear, it follows from \eqref{1.8} that
	\begin{align*}
	\norm{I_{h}(y)-I_{h}(y_{h})} \leq C (h\norm{\eta}_{1}+\norm{\eta}+h\norm{\theta}_{1}+\norm{\theta}).
	\end{align*}
	Then, by Lemmas \ref{L4.1} - \ref{L4.4}, we obtain the desired estimate.
\end{proof}
\section{Backward Euler Method}
In this section, we consider the backward Euler method for the time discretization of the semi-discrete scheme \eqref{3.1}. The error analysis of the state variable is established for a fully discrete finite element scheme. Let $ t_{n}= nk ,\medspace  \ n=0, 1, \ldots, M$ be the temporal grid points in the time interval $ [0,T],$ where $ k=\frac{T}{M} $ is the time step size. For a smooth function $ \phi $ defined on $ [0,T], $ we denote by $ \phi^{n}= \phi(t_{n})$ the value of \(\phi\) at the discrete time point \(t_{n}\).  $\bar{\partial_{t}}\phi^{n}= \frac{\phi^{n}- \phi^{n-1}}{k} $ is an expression of the time derivative.

Using the backward Euler approach on \eqref{3.1} produces a sequence of functions $ \{Y^{n}\}_{n\geq 1}\in S_{h}$ such that
\begin{align}
	\label{5.1}
	&(\bar{\partial_{t}}Y^{n}, \chi)+ \nu (Y_{x}^{n}, \chi_{x}) - \gamma(Y^{n}, \chi)+\delta ((Y^{n})^{3}, \chi )= -\mu (I_{h}(Y^{n}), \chi), \quad \text{for all} \ \ \chi\in S_{h}, \\&
\nonumber	Y^{0}=y_{0h}.
\end{align}
In the following lemma, we discuss {\it{a priori}} bounds  for the solution $ \{Y^{n}\}_{n\geq 1} $ of \eqref{5.1}.
\begin{lemma}
	\label{L5.1}
	Select \(k_{0} > 0 \) so that for \(0 < k \leq k_{0} \) we have
	\[
	e^{\alpha k}< 1+ \frac{k\beta_{3}}{2},
	\]
	where \(0\leq\beta_{3}=\mu-2\nu-2\gamma\).
	Then the discrete solution \( \{Y^{n}\}_{n\geq 1} \) of \eqref{5.1} is bounded by
	\begin{align*}
	\norm{Y^{M}}^{2} &+ k\beta_{4}e^{-\alpha (k+2 t_{M})} \sum_{n=1}^{M}\norm{\hat{Y}^{n}}_{1}^{2} \\
	&+ 2k\delta e^{-2\alpha t_{M}} \sum_{n=1}^{M} e^{-\alpha (2t_{n}+ k)} \norm{\hat{Y}^{n}}_{L^{4}}^{4} \leq e^{-2\alpha t_{M}}\norm{Y^{0}}^{2},
	\end{align*}
	where \(\beta_{4}=\min\{\beta e^{-\alpha k},\ \beta_{3} e^{-\alpha k}-2\frac{(1-e^{-\alpha k})}{k} \}>0\).
\end{lemma}
\begin{proof}
	Set \( \chi = Y^{n} \) in \eqref{5.1} to obtain
\begin{align}
\label{5.2}
(\bar{\partial_{t}}Y^{n}, Y^{n}) + \nu \norm{Y^{n}_{x}}^{2} - \gamma \norm{Y^{n}}^{2} + \delta \norm{Y^{n}}_{L^{4}}^{4} = -\mu (I_{h}(Y^{n}), Y^{n}).
\end{align}

A use of the Cauchy-Schwarz and Young's inequalities on the right-hand side of \eqref{5.2}, we obtain from \eqref{5.2} with \eqref{2.2}
\[
(\bar{\partial_{t}}Y^{n}, Y^{n}) + \frac{\beta}{2}\norm{Y_{x}^{n}}^{2}+\frac{\beta_{3}}{2}\norm{Y^{n}}^{2} + \delta \norm{Y^{n}}_{L^{4}}^{4} \leq 0,
\]
where 
\(
 \beta =2\nu-\mu c_{p}^{2}h^{2}, \  \text{and} \quad \beta_{3}=\mu - 2\nu - 2\gamma \ \text{ are positive constants}.
\)

Multiplying by \( e^{2\alpha t_{n}} \) and setting 
\(
\hat{Y}^{n} = e^{\alpha t_{n}}Y^{n},
\)
in the above inequality gives
\[
\left(e^{\alpha t_{n}}\bar{\partial_{t}}Y^{n}, \hat{Y}^{n}\right) + \frac{\beta}{2}\norm{\hat{Y_{x}}^{n}}^{2}+\frac{\beta_{3}}{2}\norm{\hat{Y}^{n}}^{2} + \delta e^{-2\alpha t_{n}} \norm{\hat{Y}^{n}}_{L^{4}}^{4} \leq 0.
\]

Since
\[
e^{\alpha t_{n}}\bar{\partial_{t}}Y^{n} = e^{\alpha k}\bar{\partial_{t}}\hat{Y}^{n} - \frac{e^{\alpha k} - 1}{k} \hat{Y}^{n},
\]
and
\[
(\bar{\partial_{t}}\hat{Y}^{n}, \hat{Y}^{n}) = \frac{1}{2}\bar{\partial_{t}} \norm{\hat{Y}^{n}}^{2} + \frac{k}{2}\norm{\bar{\partial_{t}}\hat{Y}^{n}}^{2},
\]
we obtain, after multiplying by \( 2e^{-\alpha k} \),
\[
\bar{\partial_{t}} \norm{\hat{Y}^{n}}^{2} + k\norm{\bar{\partial_{t}}\hat{Y}^{n}}^{2} -2 \frac{(1-e^{-\alpha k})}{k}\norm{\hat{Y}^{n}}^{2} + \beta e^{-\alpha k}\norm{\hat{Y_{x}}^{n}}^{2}+\beta_{3}e^{-\alpha k}\norm{\hat{Y}^{n}}^{2} + 2\delta e^{-\alpha(2t_{n}+k)} \norm{\hat{Y}^{n}}_{L^{4}}^{4} \leq 0.
\]

With \( 0 < 2\alpha \leq \beta_{3} \), select \( k_{0} > 0 \) such that for \( 0 < k \leq k_{0} \),
\[
e^{\alpha k}< 1+ \frac{k\beta_{3}}{2}.
\]
Consequently, we have
\[
\norm{\hat{Y}^{n}}^{2} + k e^{-\alpha k}\beta_{4} \norm{\hat{Y}^{n}}_{1}^{2} + 2k\delta e^{-\alpha(2t_{n}+k)} \norm{\hat{Y}^{n}}_{L^{4}}^{4} \leq \norm{\hat{Y}^{n-1}}^{2},
\]
where \(\beta_{4}=\min\{\beta e^{-\alpha k},\  \beta_{3} e^{-\alpha k}-2\frac{(1-e^{-\alpha k})}{k} \}>0\).

Summing from \( n=1 \) to \( n=M \), we observe that
\[
\norm{\hat{Y}^{M}}^{2} + k\beta_{4}e^{-\alpha k} \sum_{n=1}^{M} \norm{\hat{Y}^{n}}_{1}^{2} + 2k\delta \sum_{n=1}^{M} e^{-\alpha(2t_{n}+k)} \norm{\hat{Y}^{n}}_{L^{4}}^{4} \leq \norm{Y^{0}}^{2}.
\]

Multiplying the above inequality by \( e^{-2\alpha t_{M}} \) completes the proof.
\end{proof}

\subsection{Existence and Uniqueness of Discrete Solution.}
In this subsection, we determine the existence and uniqueness of the fully discrete scheme \eqref{5.1}.
\begin{theorem}[Brouwer's Fixed Point Theorem \cite{kesavan2023functional}]
	\label{B}
	Let \( X \) be a finite-dimensional inner product space. Suppose that \( f: X \rightarrow X \) is continuous and there exists \( \gamma > 0 \) such that
	\[
	(f(v), v) \geq 0 \quad \text{for all } v \in X \text{ with } \| v \| = \gamma.
	\]
	Then, there exists \( v_{1} \in X \) such that \( f(v_{1}) = 0 \) and \( \| v_{1} \| \leq \gamma \).
\end{theorem}

\begin{theorem}
	Given \( Y^{n-1} \), the solution of the discrete scheme \eqref{5.1} exists for sufficiently small \( k \) such that
	\[
	k\Big(\gamma+\frac{\mu c_{p}^{2}h^{2}}{2}\Big) < 1.
	\]
\end{theorem}

\begin{proof}
	For a fixed \( n \), construct a function \( g: S_{h} \rightarrow S_{h} \) given \( Y^{n-1} \) by
	\begin{align}
	\label{5.4}
	(g(V), \chi) = (V - Y^{n-1}, \chi) + k\nu (V_{x}, \chi_{x}) - k\gamma (V, \chi) + k\delta (V^{3}, \chi) + k\mu (I_{h}(V), \chi),
	\end{align}
	for all \( \chi \in S_{h} \). Since \( g \) is a polynomial function of \( V \), it is continuous.
	
	Setting \( \chi = V \) in \eqref{5.4} and using the Cauchy–Schwarz and Young's inequalities, we have from \eqref{1.8}
	\begin{align*}
	(g(V), V) &\geq \| V \|^{2} - \| Y^{n-1} \| \, \| V \| + k(\nu-\frac{\mu c_{p}^{2}h^{2}}{2}) \| V_{x} \|^{2} - k\gamma \| V \|^{2} + k\delta \| V \|_{L^{4}}^{4}- \frac{k\mu c_{p}^{2}h^{2}}{2}  \| V \|^{2}+\frac{k\mu}{2}  \| V \|^{2}.
	\end{align*}
	Therefore, we can rewrite the above inequality as
	\begin{align*}
	(g(V), V) \geq \| V \| \left( \| V \|\Big(1 - k\gamma-\frac{k\mu c_{p}^{2}h^{2}}{2}\Big) - \| Y^{n-1} \| \right) + \frac{k\beta}{2} \| V_{x} \|^{2} + k\delta \| V \|_{L^{4}}^{4}+\frac{k\mu}{2}  \| V \|^{2},
	\end{align*}
	where \(0\leq\beta =2\nu-\mu c_{p}^{2}h^{2}\).
	
	Hence, for
	\(
	\| V \| \geq \frac{\| Y^{n-1} \|}{ 1 - k\Big(\gamma+\frac{\mu c_{p}^{2}h^{2}}{2}\Big)},
	\)
	and provided that \( k\Big(\gamma+\frac{\mu c_{p}^{2}h^{2}}{2}\Big) < 1 \), we have \( (g(V), V) \geq 0 \). Therefore, by Theorem \ref{B}, there exists a \( v_{1} \in S_{h} \) such that \( g(v_{1}) = 0 \). This completes the proof.
\end{proof}

In the next theorem, we prove the uniqueness of the fully discrete scheme.
\begin{theorem}
A unique solution to the discrete problem \eqref{5.1} is admitted.
\end{theorem}

\begin{proof}
	Set 
	\[
	Z^{n} = Y_{1}^{n} - Y_{2}^{n},
	\]
	where \( Y_{1}^{n} \) and \( Y_{2}^{n} \) are two solutions of \eqref{5.1}. Then, \( Z^{n} \) satisfies
	\begin{align*}
		(\bar{\partial_{t}}Z^{n}, \chi) + \nu (Z_{x}^{n}, \chi_{x}) - \gamma (Z^{n}, \chi) + \delta \left((Y_{1}^{n})^{3} - (Y_{2}^{n})^{3}, \chi\right) = -\mu \left(I_{h}(Y_{1}^{n}) - I_{h}(Y_{2}^{n}), \chi\right)
	\end{align*}
	for all \( \chi \in S_{h} \). 

	Substituting \( \chi = Z^{n} \) into the above equation, we obtain
	\begin{align}
	\label{5.5}
		(\bar{\partial_{t}}Z^{n}, Z^{n}) + \nu \| Z_{x}^{n} \|^{2} - \gamma \| Z^{n} \|^{2} + \delta \left((Y_{1}^{n})^{3} - (Y_{2}^{n})^{3}, Z^{n}\right) = -\mu \left(I_{h}(Y_{1}^{n}) - I_{h}(Y_{2}^{n}), Z^{n}\right).
	\end{align}

	The term 
	\(
	\left((Y_{1}^{n})^{3} - (Y_{2}^{n})^{3}, Z^{n}\right)
	\)
	can be written as
	\[
	\left((Y_{1}^{n})^{3} - (Y_{2}^{n})^{3}, Z^{n}\right) = \left(Z^{n}\left((Y_{1}^{n})^{2} + Y_{1}^{n}Y_{2}^{n} + (Y_{2}^{n})^{2}\right), Z^{n}\right),
	\]
	which is nonnegative. Therefore, using \eqref{5.5},  Cauchy--Schwarz and the Young's inequalities, we obtain from \eqref{1.8}
	\[
	(\bar{\partial_{t}}Z^{n}, Z^{n}) + \nu \| Z_{x}^{n} \|^{2} - \gamma \| Z^{n} \|^{2} \leq \frac{\mu c_{p}^{2}h^{2}}{2} \| Z_{x}^{n} \|^{2}+\frac{\mu c_{p}^{2}h^{2}}{2} \| Z^{n} \|^{2}-\frac{\mu}{2}\| Z^{n} \|^{2}.
	\]
	Note that
	\[
	(\bar{\partial_{t}}Z^{n}, Z^{n}) = \frac{1}{2} \bar{\partial_{t}} \| Z^{n} \|^{2} + \frac{k}{2} \| \bar{\partial_{t}} Z^{n} \|^{2}.
	\]
	Using \eqref{2.2}, we have 
	\[
	\frac{1}{2} \bar{\partial_{t}} \| Z^{n} \|^{2} + \frac{\beta}{2} \| Z_{x}^{n} \|^{2} \leq 0,
	\]
	where \(0\leq \beta = 2\nu-\mu c_{p}^{2}h^{2}\).
	
	 From the above inequality, we observe that
	\[
	\| Z^{n} \|^{2} \leq  \| Z^{0} \|^{2}.
	\]
	 Since \( Z^{0} = Y_{1}^{0} - Y_{2}^{0} = 0 \), it follows that \( Z^{n} = 0 \) for all \( n \geq 0 \). This proves the uniqueness of the solution \( \{Y^{n}\}_{n\geq 1} \) to \eqref{5.1}.
\end{proof}
\subsection{Error Analysis of Backward Euler Method}
For the state variable and control input, the error analysis of a fully discrete scheme is examined in this subsection.

We determine the error in the following way
\[
e^{n} := y(t_{n}) - Y^{n} = \bigl( y(t_{n}) - \tilde{y}(t_{n}) \bigr) - \bigl( Y^{n} - \tilde{y}(t_{n}) \bigr) =: \eta^{n} - \theta^{n},
\]
where \( Y^{n} \)  and \( y(t_{n}) \) are the solution of  \eqref{5.1} and \eqref{1.9} respectively, at \( t = t_{n} \); here, \( \eta^{n} = \eta(t_{n}) \) is defined in (\ref{aux}) at \( t = t_{n} \). We also denote \(y^{n} = y(t_{n})\).

Using \eqref{1.9} and (\ref{aux}) at \( t = t_{n} \) and then subtracting the resulting equation from \eqref{3.1}, we obtain
\begin{align}
\label{5.7}
(\bar{\partial_{t}}\theta^{n}, \chi) + \nu (\theta_{x}^{n}, \chi_{x}) - \gamma (\theta^{n}, \chi)&+ \delta ((\theta^{n})^{3}, \chi) + 3\delta ((\theta^{n})^{2}Y^{n}, \chi) \nonumber\\[1mm]
&= (\bar{\partial_{t}}\eta^{n}, \chi) - \gamma (\eta^{n}, \chi) + \delta ((\eta^{n})^{3}, \chi)
 + 3\delta ((\eta^{n})^{2}y^{n}, \chi) - 3\delta (\eta^{n}(y^{n})^{2}, \chi) \nonumber\\[1mm]
&\quad + 3\delta (\theta^{n}(Y^{n})^{2}, \chi) 
 + \mu (I_{h}(y^{n}) - I_{h}(Y^{n}), \chi) 
 + (y_{t}(t_{n}) - \bar{\partial_{t}}y^{n}, \chi).
\end{align}
 \begin{theorem}
 	\label{th4.4}
 	Assume that \( k_{0} > 0 \) such that for \( 0 < k \leq k_{0} \)
 	\[
 	e^{\alpha k}< 1+ \frac{k\beta_{2}}{2},
 	\]
 	holds, where \( \beta_{2}=\mu-3\nu-2\gamma \geq 0\) and \(0 < 2\alpha \leq \beta_{2}\).
 	
 	Then, there exists a positive constant \( C = C(\nu, \gamma, \delta, \mu, c_{p})  \) independent of \( h \) and \( k \) such that
 	\begin{align*} 
 	\norm{\theta^{M}}^{2} &+ \beta_{5} k e^{-2\alpha t_{M}} \sum_{n=1}^{M} \norm{\hat{\theta}^{n}}_{1}^{2} + k\delta e^{-2\alpha t_{M}} \sum_{n=1}^{M} e^{-\alpha(2t_{n}+k)} \norm{\hat{\theta}^{n}}_{L^{4}}^{4} \\
 	&\leq e^{-2\alpha t_{M}} C (h^{4} + k^{2}),
 	\end{align*}
 	where 
 	\(
 	\beta= 2\nu-\mu c_{p}^{2}h^{2},\quad \beta_{5} =  \min\left\{\beta e^{-\alpha k},\, e^{-\alpha k} \beta_{2}-2\frac{(1-e^{-\alpha k})}{k} \right\}>0,
 	\)
 	and 
 	\(
 	\hat{\theta}^{n} = e^{\alpha t_{n}} \theta^{n}.
 	\)
 \end{theorem}
 \begin{proof}
 	Set \( \chi = \theta^{n} \) in \eqref{5.7} to obtain
 	\begin{align}
 	\label{5.8}
 	(\bar{\partial_{t}}\theta^{n}, \theta^{n}) &+ \nu \norm{\theta_{x}^{n}}^{2} - \gamma \norm{\theta^{n}}^{2} + \delta \norm{\theta^{n}}_{L^{4}}^{4} \nonumber\\[1mm]
 	&= (\bar{\partial_{t}}\eta^{n}, \theta^{n}) - \gamma (\eta^{n}, \theta^{n}) + \delta ((\eta^{n})^{3}, \theta^{n})  + 3\delta ((\eta^{n})^{2}y^{n}, \theta^{n}) - 3\delta (\eta^{n}(y^{n})^{2}, \theta^{n}) \nonumber\\[1mm]
 	&\quad + 3\delta (\theta^{n}(Y^{n})^{2}, \theta^{n}) - 3\delta ((\theta^{n})^{2}Y^{n}, \theta^{n})
 	 + \mu \left( I_{h}(y^{n}) - I_{h}(Y^{n}), \theta^{n} \right) + \left( y_{t}(t_{n}) - \bar{\partial_{t}}y^{n}, \theta^{n} \right), \nonumber\\[1mm]
 	&=: \sum_{i=1}^{4} T_{i}(\theta^{n}).
 	\end{align}	
  On the right-hand side, the first term is estimated by
 	\begin{align*}
 	T_{1}(\theta^{n}) &= \left( (\bar{\partial_{t}}\eta^{n}, \theta^{n}) - \gamma (\eta^{n}, \theta^{n}) \right) + \delta ((\eta^{n})^{3}, \theta^{n}), \\[1mm]
 	&\leq \frac{\nu}{8}\norm{\theta^{n}}^{2} + C \norm{\bar{\partial_{t}}\eta^{n}}^{2} + C(\gamma, \nu) \norm{\eta^{n}}^{2} + C \norm{\eta^{n}}_{L^{4}}^{4} + \frac{\delta}{4}\norm{\theta^{n}}_{L^{4}}^{4}.
 	\end{align*}
 	Using Young's inequality with \( Y^{n} = \tilde{y}^{n} + \theta^{n} \), the second term is estimated by
 	\begin{align*}
 	T_{2}(\theta^{n}) &= 3\delta \bigl( (\theta^{n})^{2}\tilde{y}^{n}, \theta^{n} \bigr) + 3\delta \bigl( \theta^{n} (\tilde{y}^{n})^{2}, \theta^{n} \bigr) - 3\delta \bigl( \eta^{n}(y^{n})^{2}, \theta^{n} \bigr) + 3\delta \bigl( (\eta^{n})^{2}y^{n}, \theta^{n} \bigr), \\[1mm]
 	&\leq \frac{\nu}{8}\norm{\theta^{n}}^{2} + \frac{\delta}{4}\norm{\theta^{n}}_{L^{4}}^{4} + C(\delta)\norm{\tilde{y}^{n}}_{\infty}^{2}\norm{\theta^{n}}^{2} + C(\delta)\norm{y^{n}}_{\infty}^{4}\norm{\eta^{n}}^{2} + C(\delta, \nu)\norm{y^{n}}_{\infty}^{2}\norm{\eta^{n}}_{L^{4}}^{4}.
 	\end{align*}	
 	The third term, as in the semi-discrete case, is bounded by
 	\begin{align*}
 	T_{3}(\theta^{n}) &= \mu \left( I_{h}(y^{n}) - I_{h}(Y^{n}), \theta^{n} \right),
 	 \\&\leq C(\mu, \nu, c_{p})h^{2}\norm{\eta^{n}}_{1}^{2}+C(\mu, \nu)\norm{\eta^{n}}^{2}+\frac{\nu}{8}\norm{\theta^{n}}^{2}+\frac{\mu c_{p}^{2}h^{2}}{2}\norm{\theta^{n}_{x}}^{2}+\frac{\mu c_{p}^{2}h^{2}}{2}\norm{\theta^{n}}^{2}-\frac{\mu}{2}\norm{\theta^{n}}^{2}.
 	 \end{align*}	
 	Finally, with the help of Young's inequality, the fourth term yields
 	\begin{align*}
 	T_{4}(\theta^{n}) &= \bigl( y_{t}(t_{n}) - \bar{\partial_{t}}y^{n}, \theta^{n} \bigr)\leq \frac{\nu}{8}\norm{\theta^{n}}^{2} + \frac{2}{\nu}\norm{y_{t}(t_{n}) - \bar{\partial_{t}}y^{n}}^{2}.
 	\end{align*}
 	
 	Substituting these estimates into \eqref{5.8} and using the facts that \( H^{1} \xhookrightarrow{} L^{\infty} \) and 
 	\[
 	\norm{\tilde{y}^{n}}_{\infty}^{2} \leq C \norm{y^{n}}_{1}^{2},
 	\]
 	after multiplying by \( e^{2\alpha t_{n}} \) and setting 
 	\(
 	\hat{\theta}^{n} = e^{\alpha t_{n}}\theta^{n},
 	\)
 	we arrive at
 	\begin{align}
 	\label{5.9}
 	e^{\alpha t_{n}}(\bar{\partial_{t}}\theta^{n}, \hat{\theta}^{n}) + \frac{\beta}{2} \norm{\hat{\theta}_{x}^{n}}^{2} + \frac{\delta e^{-2\alpha t_{n}}}{2}\norm{\hat{\theta}^{n}}_{L^{4}}^{4} 
 	&\leq \left(\frac{\nu}{2} + \frac{\mu c_{p}^{2}h^{2}}{2} + \gamma-\frac{\mu}{2}\right) \norm{\hat{\theta}^{n}}^{2} \nonumber\\[1mm]
 	&\quad + C e^{2\alpha t_{n}}\Bigl( \norm{\bar{\partial_{t}}\eta^{n}}^{2} + \norm{\eta^{n}}^{2} + \norm{\eta^{n}}_{L^{4}}^{4}+h^{2}\norm{\eta^{n}}_{1}^{2} \Bigr) \nonumber\\[1mm]
 	&\quad + C \norm{y^{n}}_{1}^{2}\norm{\hat{\theta}^{n}}^{2} + \frac{2e^{2\alpha t_{n}}}{\nu}\norm{y_{t}(t_{n}) - \bar{\partial_{t}}y^{n}}^{2},
 	\end{align}
 	where \(\beta = 2\nu- \mu c_{p}^{2}h^{2}\) is a nonnegative constant.
 	
 		Also
 	\begin{align}
 	\label{5.10}
 	e^{\alpha t_{n}}\bar{\partial_{t}}\theta^{n} = e^{\alpha k}\bar{\partial_{t}}\hat{\theta}^{n} - \left( \frac{e^{\alpha k} - 1}{k} \right)\hat{\theta}^{n}.
 	\end{align}
 	
 	Using \eqref{5.10} and \eqref{2.2} in \eqref{5.9} gives
 	\begin{align}
 	\label{5.11}
 	e^{\alpha k}(\bar{\partial_{t}}\hat{\theta}^{n}, \hat{\theta}^{n}) &- \left( \frac{e^{\alpha k} - 1}{k} \right)\norm{\hat{\theta}^{n}}^{2} + \frac{\beta}{2} \norm{\hat{\theta}_{x}^{n}}^{2} + \frac{\delta e^{-2\alpha t_{n}}}{2}\norm{\hat{\theta}^{n}}_{L^{4}}^{4} \nonumber\\[1mm]
 	&\leq \left(\frac{3\nu}{2} + \gamma - \frac{\mu}{2} \right) \norm{\hat{\theta}^{n}}^{2} \nonumber\\[1mm]
 	&\quad + C e^{2\alpha t_{n}}\Bigl( \norm{\bar{\partial_{t}}\eta^{n}}^{2} + \norm{\eta^{n}}^{2} + \norm{\eta^{n}}_{L^{4}}^{4} +h^{2}\norm{\eta^{n}}_{1}^{2}\Bigr) \nonumber\\[1mm]
 	&\quad + C \norm{y^{n}}_{1}^{2}\norm{\hat{\theta}^{n}}^{2} + \frac{2e^{2\alpha t_{n}}}{\nu}\norm{y_{t}(t_{n}) - \bar{\partial_{t}}y^{n}}^{2}.
 	\end{align}	
 	Since we can write 
 	\begin{align*}
 	(\bar{\partial_{t}}\hat{\theta}^{n}, \hat{\theta}^{n}) = \frac{1}{2}\bar{\partial_{t}}\norm{\hat{\theta}^{n}}^{2} + \frac{k}{2}\norm{\bar{\partial_{t}}\hat{\theta}^{n}}^{2},
 	\end{align*}
 	it follows from \eqref{5.11} that
 	\begin{align}
 	\label{5.12}
 	\frac{e^{\alpha k}}{2}\bar{\partial_{t}}\norm{\hat{\theta}^{n}}^{2} + \frac{k}{2}\norm{\bar{\partial_{t}}\hat{\theta}^{n}}^{2} &- \left( \frac{e^{\alpha k} - 1}{k} \right)\norm{\hat{\theta}^{n}}^{2} + \frac{\beta}{2} \norm{\hat{\theta}_{x}^{n}}^{2} + \frac{\delta e^{-2\alpha t_{n}}}{2}\norm{\hat{\theta}^{n}}_{L^{4}}^{4} \nonumber\\[1mm]
 	&\leq \left(\frac{3\nu}{2} + \gamma - \frac{\mu}{2} \right)\norm{\hat{\theta}^{n}}^{2} \nonumber\\[1mm]
 	&\quad + C e^{2\alpha t_{n}}\Bigl( \norm{\bar{\partial_{t}}\eta^{n}}^{2} + \norm{\eta^{n}}^{2} + \norm{\eta^{n}}_{L^{4}}^{4}+h^{2}\norm{\eta^{n}}_{1}^{2} \Bigr) \nonumber\\[1mm]
 	&\quad + C \norm{y^{n}}_{1}^{2}\norm{\hat{\theta}^{n}}^{2} + \frac{2e^{2\alpha t_{n}}}{\nu}\norm{y_{t}(t_{n}) - \bar{\partial_{t}}y^{n}}^{2}.
 	\end{align}
 	
 	With \(0 < 2\alpha \leq \beta_{2}\), select \( k_{0} > 0 \) such that for \( 0 < k \leq k_{0} \)
 	\[
 	e^{\alpha k}< 1+ \frac{k\beta_{2}}{2},
 	\]
 	where \(\beta_{2}= \mu-3\nu -2\gamma\) is a positive constant.
 	
 	 Multiplying \eqref{5.12} by \( 2e^{-\alpha k} \) yields
 	\begin{align}
 	\label{5.13}
 	\bar{\partial_{t}}\norm{\hat{\theta}^{n}}^{2} + \beta_{5} \norm{\hat{\theta}^{n}}_{1}^{2} &+ \delta e^{-\alpha(2t_{n}+k)} \norm{\hat{\theta}^{n}}_{L^{4}}^{4} \nonumber\\[1mm]& \leq e^{-\alpha k} C\norm{y^{n}}_{1}^{2}\norm{\hat{\theta}^{n}}^{2}  + Ce^{\alpha(2t_{n}-k)}\left(\norm{\bar{\partial_{t}}\eta^{n}}^{2} + \norm{\eta^{n}}^{2} + \norm{\eta^{n}}_{L^{4}}^{4}+h^{2}\norm{\eta^{n}}_{1}^{2}\right) \nonumber\\[1mm]& 
 	\quad+ \frac{4e^{\alpha(2t_{n}-k)}}{\nu}\norm{y_{t}(t_{n}) - \bar{\partial_{t}}y^{n}}^{2},
 	\end{align}
 	where
 	\(
 	\beta_{5} =  \min\left\{\beta e^{-\alpha k},\, e^{-\alpha k} \beta_{2}-2\frac{(1-e^{-\alpha k})}{k} \right\}>0.
 	\)
 	
 	In the interval \( (t_{n-1}, t_{n}) \), we utilize the Taylor series expansion of \( y \) about the point \( t_{n} \) to get
 	\begin{align}
 	\label{5.14}
 	\norm{y_{t}(t_{n}) - \bar{\partial_{t}}y^{n}}^{2} \leq \frac{k}{3} \int_{t_{n-1}}^{t_{n}} \norm{y_{tt}(s)}^{2} ds.
 	\end{align}
 	
 	Summing \eqref{5.13} from \( n=1 \) to \( n=M \) and multiplying \eqref{5.14} by \( k \) yields
 	\begin{align}
 	\label{5.15}
 	\norm{\hat{\theta}^{M}}^{2} + \beta_{5} k \sum_{n=1}^{M}\norm{\hat{\theta}^{n}}_{1}^{2} &+ k\delta \sum_{n=1}^{M}e^{-\alpha(2t_{n}+k)}\norm{\hat{\theta}^{n}}_{L^{4}}^{4}  \nonumber\\[1mm]&\leq\; k e^{-\alpha k} \sum_{n=1}^{M}  C \norm{y^{n}}_{1}^{2} \norm{\hat{\theta}^{n}}^{2}+C(\nu)k^{2}\int_{0}^{t_{M}}\norm{y_{tt}(s)}^{2} ds \nonumber\\[1mm]
 	&\quad+ Ck \sum_{n=1}^{M} e^{\alpha(2t_{n}-k)}\left(\norm{\bar{\partial_{t}}\eta^{n}}^{2} + \norm{\eta^{n}}^{2} + \norm{\eta^{n}}_{L^{4}}^{4}+h^{2}\norm{\eta^{n}}_{1}^{2}\right).
 	\end{align}
 	
 	Note that
 	\begin{align}
 	\label{5.16}
 	\norm{\bar{\partial_{t}}\eta^{n}}^{2} = \frac{1}{k^{2}} \left( \int_{t_{n-1}}^{t_{n}} (y-\tilde{y})_{t} ds \right)^{2} \leq \frac{1}{k} \int_{t_{n-1}}^{t_{n}} \norm{\eta_{t}}^{2} ds.
 	\end{align}
 	
 	Using Lemmas \ref{L4.1} and \ref{L5.1} together with \eqref{5.16} and applying the discrete Gronwall inequality, then multiplying by \( e^{-2\alpha t_{M}} \) in the resulting inequality, we obtain
 	\begin{align*}
 	\norm{\theta^{M}}^{2} &+ \beta_{5} k e^{-2\alpha t_{M}} \sum_{n=1}^{M}\norm{\hat{\theta}^{n}}_{1}^{2} + k\delta e^{-2\alpha t_{M}} \sum_{n=1}^{M} e^{-\alpha(2t_{n}+k)}\norm{\hat{\theta}^{n}}_{L^{4}}^{4} \\
 	&\leq Ce^{-2\alpha t_{M}} \Biggl( Ch^{4}e^{-\alpha k} \left(k\sum_{n=1}^{M} e^{2\alpha t_{n}} \right)\Bigl( \norm{y(t_{n})}_{2}^{2} + \norm{y(t_{n})}_{2}^{4} +h^{2} \norm{y(t_{n})}_{2}^{2}\Bigr) \\
 	&\quad +Ch^{4} \int_{0}^{t_{M}}\norm{y_{t}(s)}_{2}^{2}ds+ C(\nu) k^{2} \int_{0}^{t_{M}} \norm{y_{tt}(s)}^{2} ds \Biggr).
 	\end{align*}
 	Note that 
 	\[ e^{-2\alpha t_{M}}\left(k\sum_{n=1}^{M} e^{2\alpha t_{n}} \right)=e^{-2\alpha t_{M}}\frac{1}{e^{2\alpha k}-1} ke^{2\alpha k}(e^{2\alpha t_{M}}-1)\leq C.\]
 	The use of Theorem \ref{th2.1} with \( \alpha = 0 \) completes the proof.
 \end{proof}
\begin{remark}
In the case of Neumann boundary conditions, using the auxiliary projection from Remark \ref{r3.2} at \( t = t_{n} \), we write \eqref{5.7} in the following form:
\begin{align*}
(\bar{\partial_{t}}\theta^{n}, \chi) + \nu (\theta_{x}^{n}, \chi_{x}) -& \gamma (\theta^{n}, \chi) + \delta ((\theta^{n})^{3}, \chi) + 3\delta ((\theta^{n})^{2}Y^{n}, \chi) \\
&= (\bar{\partial_{t}}\eta^{n}, \chi) - (\nu\lambda + \gamma)(\eta^{n}, \chi) + \delta ((\eta^{n})^{3}, \chi)
 + 3\delta ((\eta^{n})^{2}y^{n}, \chi) - 3\delta (\eta^{n}(y^{n})^{2}, \chi)\\[1mm]
&\quad + 3\delta (\theta^{n}(Y^{n})^{2}, \chi) + \mu \bigl(I_{h}(y^{n}) - I_{h}(Y^{n}), \chi\bigr)
+ \bigl(y_{t}(t_{n}) - \bar{\partial_{t}}y^{n}, \chi\bigr).
\end{align*}
The error analysis of the fully discrete scheme for Neumann boundary conditions is conducted in a manner analogous to that for mixed boundary conditions. For Dirichlet boundary conditions, the error analysis follows from Theorem \ref{th4.4}.
\end{remark}

\begin{theorem}
\label{th4.5}
Suppose that the assumptions of Theorem \ref{th4.4} hold. Then, there exists a positive constant \( C \) independent of \( h \) and \( k \) such that
\[
\norm{y^{n} - Y^{n}} \leq C \bigl( h^{2} + k \bigr).
\]
\end{theorem}

\begin{proof}
	The proof is obtained from Theorem \ref{th4.4} using the triangle inequality.
\end{proof}

In the following theorem, we analyze the error in the control input for the fully discrete setup. Although the theorem can be proved under assumption \eqref{1.8} and by using the bound on \( \|\theta_{x}^{n}\| \), proving this bound separately would be required. To simplify the error analysis for the control variable, we assume that the interpolation operator \( I_h \) is stable, as stated below.
\begin{theorem}
	\label{th4.6}
	Assume that the hypotheses of Theorem \ref{th4.4} hold. Further, suppose that the interpolation operator \( I_h \) is stable in \( L^{2} \), i.e.,
	\begin{align}
		\label{4.16}
		\|I_{h}(\phi)\| \leq c_{0} \|\phi\|,
	\end{align}
	where $c_{0}$  is a positive constant.
	Then, there exists a positive constant \( C \) independent of \( h \) and \( k \) such that
	\[
	\norm{I_{h}(y^{n}) - I_{h}(Y^{n})} \leq C (h^{2} + k).
	\]
\end{theorem}
\begin{proof}
	Since the interpolation operator is linear, it follows from \eqref{4.16} that
	\[
	\norm{I_{h}(y^{n}) - I_{h}(Y^{n})} \leq c_{0} \norm{y^{n} - Y^{n}}.
	\]
	The proof is completed by applying Theorem \ref{th4.5}.
\end{proof}

\section{Numerical Examples}
In this section, we present numerical examples to demonstrate stabilization and error decay. We analyze the finite parameter feedback control algorithm for the Chafee-Infante (CI) equation \eqref{1.4} - \eqref{1.6}. We show that the approximate solution of the fully discrete scheme \eqref{5.1} converges to zero, implying that the numerical solution of \eqref{1.4} - \eqref{1.6} approaches the steady-state solution. Additionally, for both state variable and the control input, we determine the order of convergence. 

To solve the resulting nonlinear system of equations for the unknown solution \( Y^{n+1} \) at each time step, we employ Newton’s method, using the solution from the previous time step \( Y^{n} \) as the initial guess.

\begin{example}
	We select the initial condition \( y_{0}(x) \) as \( y_{0}(x) = x(1 - x) \), with parameters \( \nu = 0.1, c_{p}=1, h=0.01, \) and \( \gamma = \delta = 9 \). 
\end{example}
For the uncontrolled problem \eqref{1.1}-\eqref{1.3}, the steady-state solution \( y^{\infty} = 0 \) is unstable, meaning the solution does not converge to the steady-state solution, while \( y^{\infty} = \pm 1 \) are stable. The problem \eqref{1.1}–\eqref{1.14} is asymptotically stable around the steady-state solution \( y^{\infty} = 0 \) in the absence of control whenever \( \nu \geq \gamma \). However, for the zero Neumann boundary case—i.e., the system \eqref{1.1}–\eqref{1.2}, \eqref{1.15}—asymptotic stability cannot be concluded without control. Furthermore, when \( \nu < \gamma \), the steady state \( y^{\infty} = 0 \) is not necessarily stable under the zero boundary conditions \eqref{1.3}–\eqref{1.15}, highlighting the need for feedback control to ensure stabilization.
 In contrast, with the application of a finite parameter feedback control, the equation becomes stable for sufficiently large values of \( \mu > 0 \). 

\begin{figure}[ht!]
		\centering
		(i) \includegraphics[width= 0.30\textwidth]{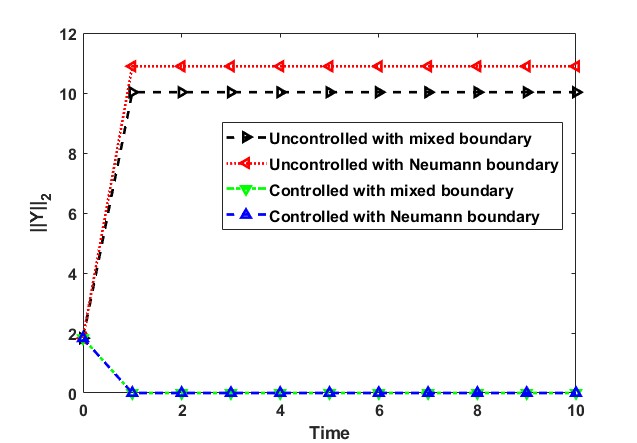}
		(ii) \includegraphics[width= 0.60\textwidth]{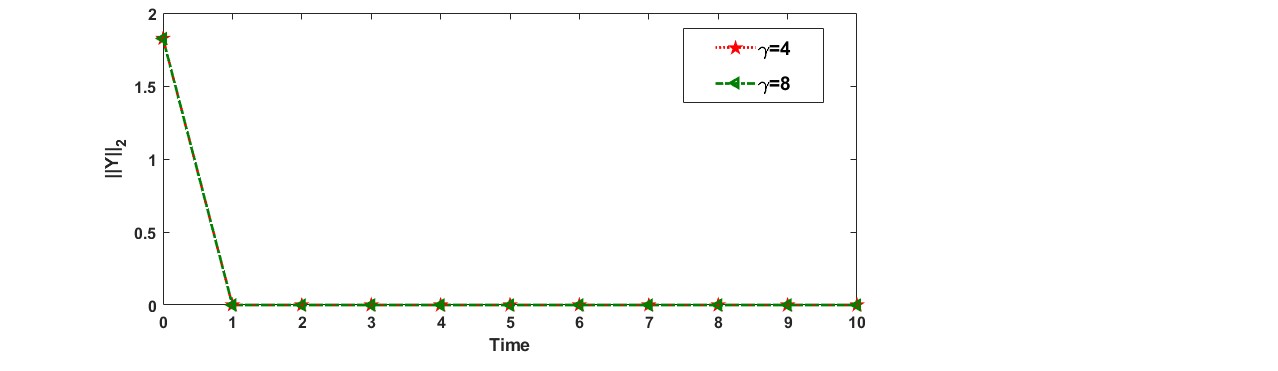}	
		(iii)\includegraphics[width= 0.25\textwidth]{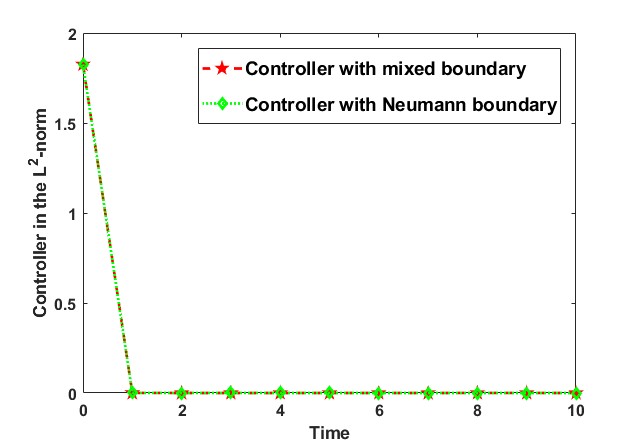}
		(iv)\includegraphics[width= 0.65\textwidth]{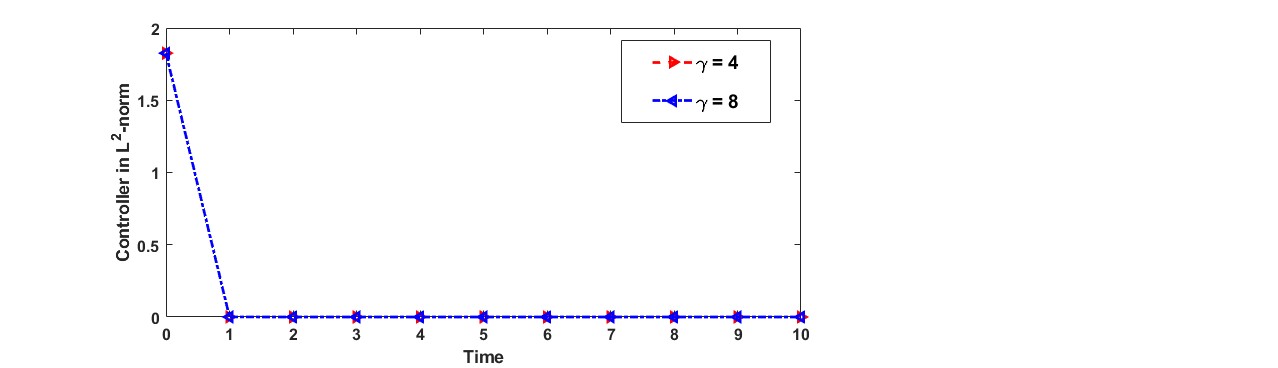}
	    (v) \includegraphics[width= 0.50\textwidth]{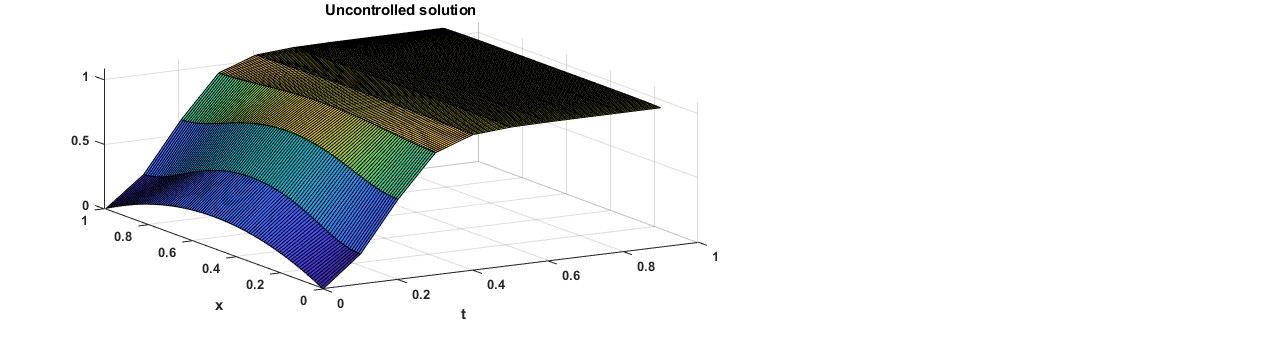}
	    (vi)\includegraphics[width= 0.40\textwidth]{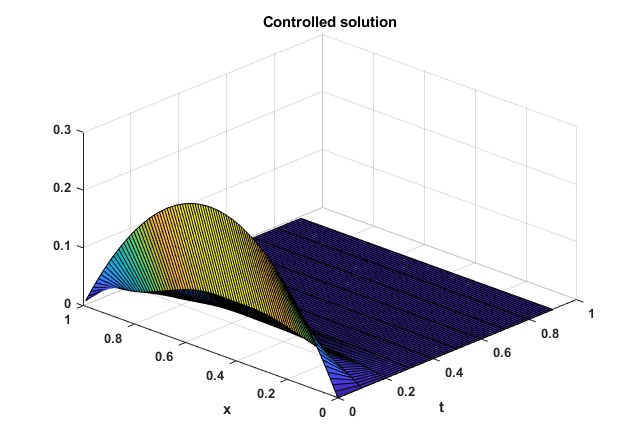}
	      
		\caption{Example 5.1: {\bf (i)} State variables in the \( L^{2} \)-norm for both mixed and Neumann boundaries with \( \gamma = 9 \), \( \mu = 20 \), and \( \nu = 0.1 \).  
{\bf (ii)} Solution of the state variable in the \( L^{2} \)-norm for different values of \( \gamma \), with fixed \( \nu = 0.1, \delta=9,\) and \( \mu = 20 \).  
{\bf (iii)} Control input for both mixed and Neumann boundaries.  
{\bf (iv)} Control input for various values of \( \gamma \).  
{\bf (v)} Uncontrolled solution with \( \gamma = 9 \), \( \mu = 20 \), and \( \nu = 0.1 \).  
{\bf (vi)} Controlled solution with \( \gamma = 9 \), \( \mu = 20 \), and \( \nu = 0.1 \).}  
\label{fig:ex1}
	\end{figure}

Here, we implement a finite parameter feedback control algorithm on nodal values to stabilize the steady-state solution of \eqref{1.1} - \eqref{1.3}.	
   In Figure \ref{fig:ex1} (i), the trajectory plot for the mixed and Neumann boundaries in the \( L^{2} \)-norm, using the same parameters, illustrates the behavior of the state variable. The behavior of the state variable for both boundaries appears to be nearly identical. Furthermore, we observe that the state variable \( Y^{n} \) in the \( L^{2} \)-norm fails to converge to a stationary solution in the absence of feedback control for both boundaries, as depicted in Figure \ref{fig:ex1} (i). 

Applying the feedback control law with a feedback parameter \( \mu = 20 \) results in the state variable approaching zero for both boundaries, as shown in Figure \ref{fig:ex1} (i). Figure \ref{fig:ex1} (ii) presents the state variable for various values of \( \gamma \) in the \( L^{2} \)-norm. Moreover, the state variable converges to a stationary solution with a decay rate \( \alpha \) satisfying \eqref{2.1}.  

Figures \ref{fig:ex1} (iii) and (iv) demonstrate the behavior of the control input for both mixed and Neumann boundaries, as well as for different values of \( \gamma \), respectively. Figures \ref{fig:ex1} (v) and (vi) present surface plots of both the uncontrolled and controlled solutions, respectively, in the case of a mixed boundary condition.

	The order of convergence for the state variable in space at time $T=1$, corresponding to various values of $h$ with a fixed $M$, is displayed in Table \ref{table:1}. The feedback control parameter is set to $\mu=20$. Since the exact solution of the CI equation is unknown, we use a reference solution as an approximation of the exact solution. Specifically, we obtain the reference solution at $h=\frac{1}{1280}$ with $M=1050$.  

From Table \ref{table:1}, we conclude that the state variable is of order $2$ in both the $L^{2}$-norm and $L^{\infty}$-norm, verifying our results in Theorem \ref{th4.5}. Moreover, we achieve second-order convergence in space  satisfying the stabilization condition \eqref{2.2}.  

Table \ref{table:2} presents the order of convergence for the state variable with respect to time, considering various values of $M$ while keeping $h$ fixed. From Table \ref{table:2}, we observe that the state variable is of order $1$, in agreement with Theorem \ref{th4.5}. 

Lastly, Table \ref{table:4} illustrates the order of convergence for the control input in the $L^{\infty}$-norm, corresponding to both space and time. These results confirm the theoretical findings in Theorem \ref{th4.6}.

	\begin{table}[H]
	\caption{ The order of convergence (O.C.) with respect to space in Example 5.1, considering varying values of $h$ and a fixed value of $M=1050$.}
	\begin{tabular}{ c c c c c }
	\toprule
		h & $L^{2}-$ Error         & O. C. & $L^{\infty}-$Error    &          O. C.    \\
	 \midrule
		$\frac{1}{10}$ &   $ 6.0910e-06$       &     $--  $     &  $9.2296e-6 $ &    $--$                 \\
	           
		$\frac{1}{20}$&   $   1.8789e-06 $     &  $ 1.69 $     &   $3.0597e-06$ &    $1.59$           \\
		          
		$\frac{1}{40}$ &       $5.3748e-07 $      & $1.81$       &    $9.2419e-07$&     $1.73$              \\
	
		$\frac{1}{80}$ &      $  1.3807e-07  $     &  $1.96$       &  $2.3895e-07$  &  $1.95$           \\
		
		$\frac{1}{160}$ &        $3.4133e-08 $    & $ 2.02 $         &  $5.9380e-08$  &      $2.01$            \\
	
		$\frac{1}{320}$ &        $ 8.0383e-09$    & $2.08$         &  $1.4084e-08$  &      $2.07$            \\
	\bottomrule
	\end{tabular}
	\label{table:1}
\end{table}

\begin{table}[H]
	\centering
	\caption{The order of convergence (O.C.) with respect to time in Example 5.1 for varying values of $M$ and a fixed value of $h$.}
	\begin{tabular}{ c c c c  c }
		\toprule
		M & $L^{\infty}-$ Error($\gamma=5, h=0.005$)        & O. C. & $L^{\infty}-$Error ($\gamma=9, h=0.005$)    &          O. C.    \\
		\midrule
		$100$ &   $ 1.478e-07$       &     $--  $     & $5.1667e-6 $ &    $--$                 \\
		
		$200$&   $   8.9536e-08 $     &  $ 0.724 $     & $3.3266e-06$ &    $0.635$           \\
		
		$400$ &       $4.9723e-08 $      & $0.847$       &  $1.8827e-06$&     $0.821$              \\
		
		$800$ &      $  2.6017e-08  $     &  $0.935$      & $9.9158e-07$  &  $0.925$           \\
		
		$1600$ &        $1.3041e-08 $    & $ 0.996 $         &$4.9829e-07$  &      $0.993$            \\
		$3200$ &        $6.2479e-09 $    & $ 1.06 $         &$2.3900e-07$  &      $1.06$            \\
		\bottomrule
		
	\end{tabular}
	\label{table:2}
\end{table}

\begin{table}[H]
	\centering
	\caption{The order of convergence (O.C.) for the controller with respect to space and time in Example 5.1 for varying values of $h$ and $M$, respectively.}
	\begin{tabular}{ c c c c c c }
		\toprule
		$ h $ & $L^{\infty}-$ Error($ M=1050$)       & O. C.& M & $L^{\infty}-$Error($\gamma=9, h=0.005$)   &          O. C.    \\
		\midrule
		$\frac{1}{10}$ &   $ 1.3728e-04$       &     $--  $     & 100& $9.8941e-05 $ &    $--$                 \\
		
		$\frac{1}{20}$&   $  3.2299e-05 $     &  $ 2.08 $     &  200& $6.1820e-05$ &    $0.679$           \\
		
		$\frac{1}{40}$ &       $8.9682e-06 $      & $1.85$  & 400&$3.4468e-05$&     $0.843$              \\
		
		$\frac{1}{80}$ &      $ 2.2890e-06 $     &  $1.97$    & 800& $1.8019e-05$  &  $0.936$           \\
		
		$\frac{1}{160}$ &        $5.6525e-07 $    & $ 2.01 $   & 1600&$9.0213e-06$  &      $0.998$            \\
		
		$\frac{1}{320}$ &        $1.3311e-07 $    & $ 2.08 $ &  3200&$4.3189e-06$  &      $1.063$            \\
		\bottomrule
		
	\end{tabular}
	\label{table:4}
\end{table}

\begin{figure}[H]
		\centering
		(i)\includegraphics[width= 0.45\textwidth]{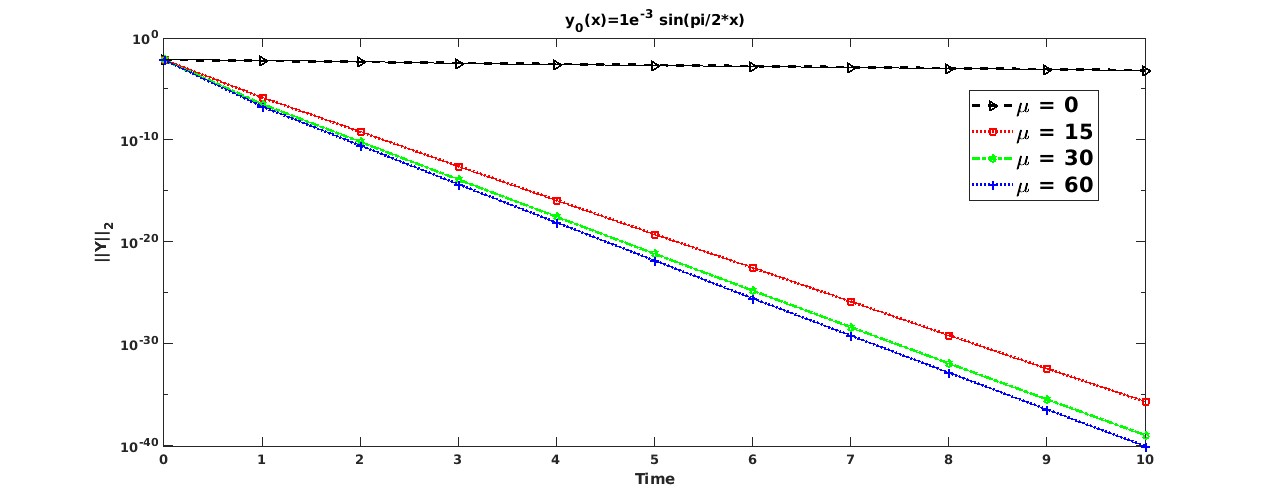}
		(ii) \includegraphics[width= 0.45\textwidth]{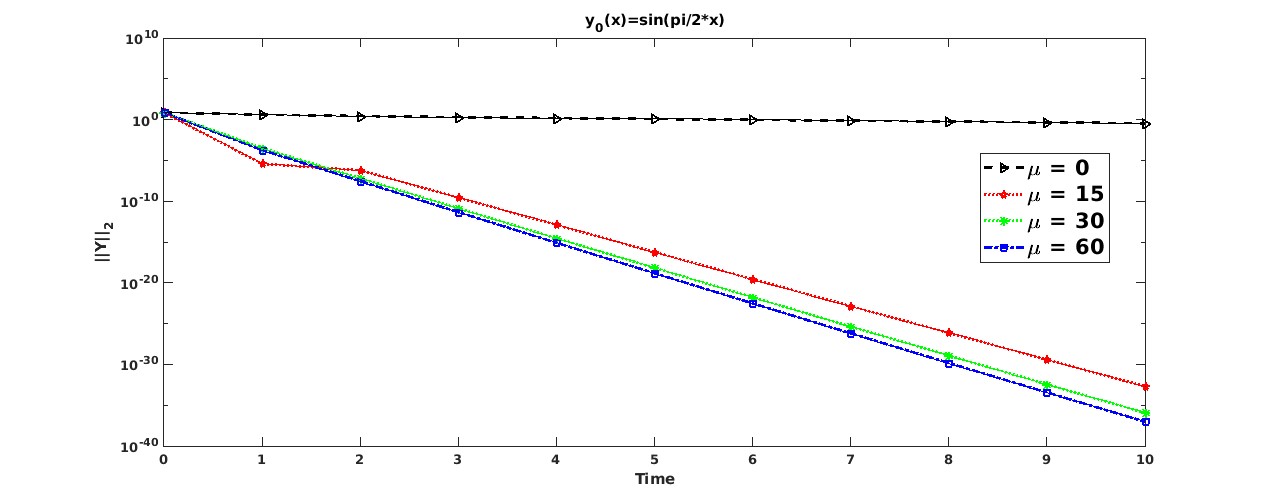}
		(iii) \includegraphics[width= 0.44\textwidth]{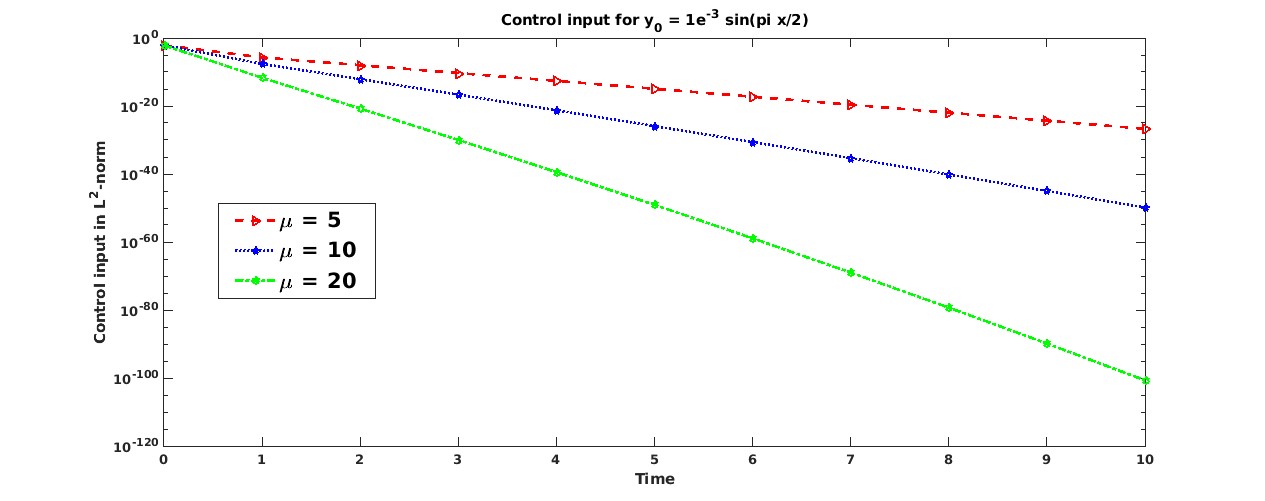}
		(iv) \includegraphics[width= 0.45\textwidth]{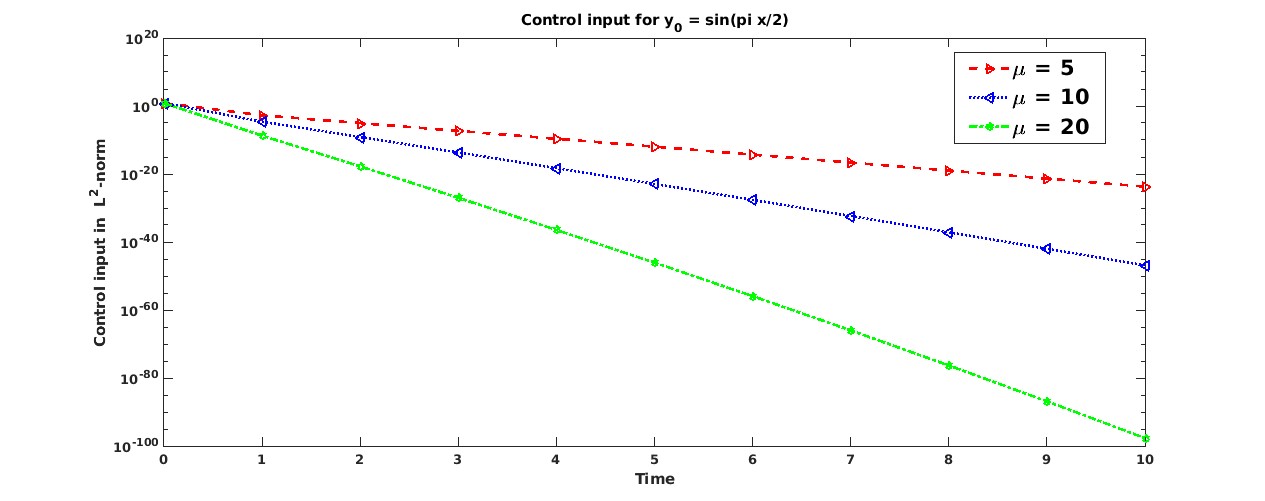}
		(v) \includegraphics[width= 0.45\textwidth]{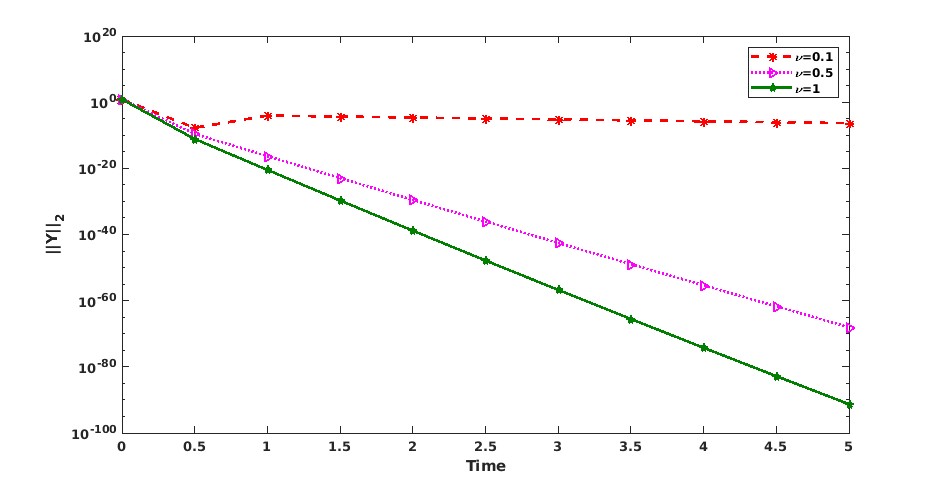}
		(vi) \includegraphics[width= 0.45\textwidth]{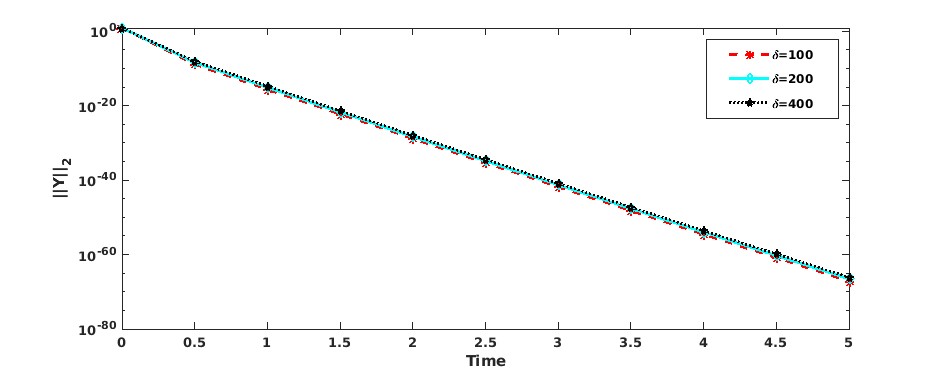}
		\caption{Example 5.2: {\bf(i)} and {\bf(ii)} Uncontrolled and controlled solutions for different values of $\mu$ with fixed parameters $\nu = 0.5$ and $\gamma = \delta = 1$, {\bf(iii)} and {\bf(iv)} Control input for different values of $\mu$ with fixed parameters $\nu = 0.5$ and $\gamma = \delta = 1$, {\bf(v)} State variable in the $L^{2}-$norm for different values of $\nu$ with fixed $\gamma = 50$ and $\mu = 120$, {\bf(vi)} State variable in the $L^{2}-$norm for different values of $\delta$ with fixed parameters $\gamma = 50$, $\nu = 0.5$, and $\mu = 120$.}\label{fig:ex2}
	\end{figure}

In the following example, we discuss the stabilization of the state variable for different values of $\mu$. Moreover, we examine the behavior of the state variable for various values of $\nu$ and the coefficient of the cubic reaction term $\delta$. We also use a finite feedback controller at the nodal values.
\begin{example}
In this example, we examine two types of initial conditions at $t=0$, namely, $ y_{01}(x)$ and $ y_{02}(x),$ where $ y_{01}(x) = \sin(\frac{\pi}{2}x) $ and $ y_{02}(x) = 1e^{-3} \sin(\frac{\pi}{2}x)$. We choose the parameter values as $\nu = 0.5$ and $\gamma = \delta = 1$.
\end{example}
Figures \ref{fig:ex2} (i) and (ii) display the behavior of uncontrolled ($\mu = 0$) and controlled solutions for different values of $\mu$ corresponding to the initial conditions $y_{02}(x)$ and $y_{01}(x)$, respectively. We observe that the behavior of the solutions remains nearly the same in the $ L^{2}-$norm for large values of $\mu$. Figures \ref{fig:ex2} (iii) and (iv) illustrate the behavior of the control inputs at nodal values for different values of $\mu$, corresponding to the initial conditions $y_{02}(x)$ and $y_{01}(x)$, respectively.

For fixed large values of $\gamma = \delta = 50$ and $\mu = 120$, given the initial condition $y_{01}(x)$, Figure \ref{fig:ex2} (v) shows how the state variables behave for different values of $\nu$. We observe that the state variable takes more time to reach a stationary solution for smaller values of $\nu$. Moreover, Figure \ref{fig:ex2} (vi) illustrates the behavior of the state variable for different values of the reaction coefficient $\delta$.

In the next example, we discuss the number of actuators corresponding to the interpolant operators. Furthermore, we compare the results for different interpolant operators over the time interval $t = [0,5]$.

\begin{figure}[!h]
		\centering
		(i)\includegraphics[width= 0.44\textwidth]{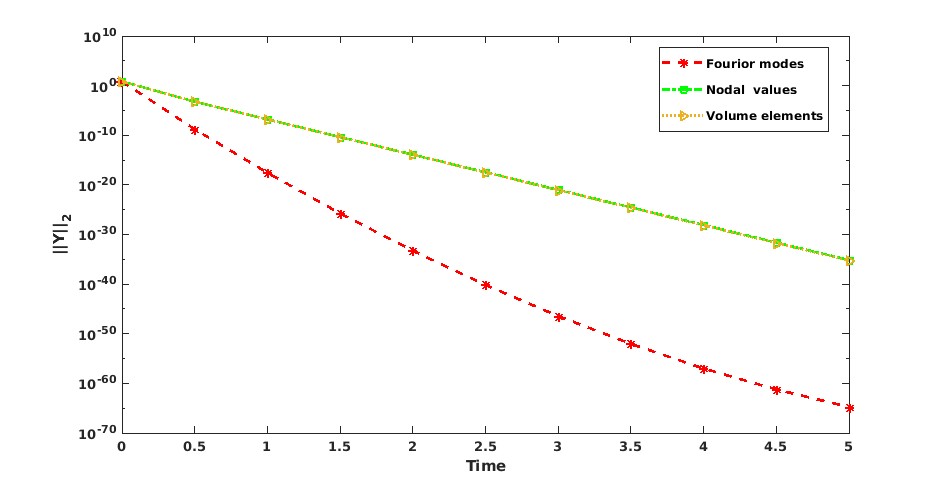}
		(ii)\includegraphics[width= 0.48\textwidth]{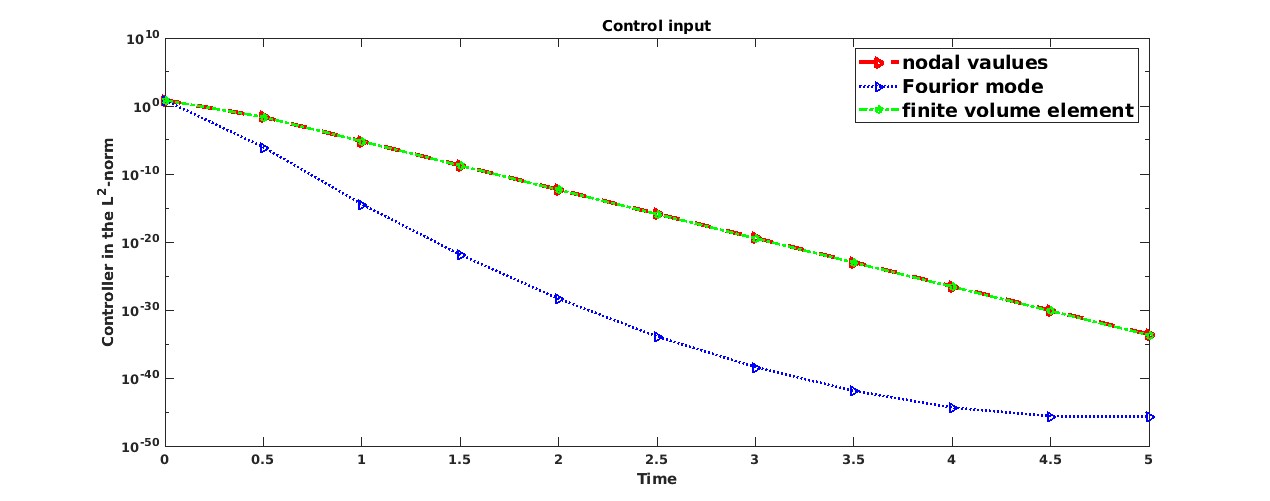}
		\caption{Example 5.3: {\bf (i)} Interpolant operators such as Fourier modes, nodal values, and finite volume elements in the $ L^{2}-$norm, {\bf (ii)} Controller for different values of finite parameter feedback control in the $ L^{2}-$norm.}
		\label{fig:ex3}
	\end{figure}

\begin{example}
	Here, we start with the initial condition \( y_{0}(x) = \cos(3\pi x) \), with fixed parameters \( \nu=1 \), \( \gamma=\delta=150 \), and \( \mu=500 \). We choose \( h=0.01 \) with a sufficiently small value of \( k \).  
\end{example}
Table \ref{table:3} presents the number of controllers corresponding to different interpolant operators, such as Fourier modes, nodal values, and volume elements. Additionally, from Figure \ref{fig:ex3} (i), we observe that the state variable in the \( L^{2} \)-norm behaves similarly for Fourier modes, nodal values, and volume elements. Figure \ref{fig:ex3} (ii) illustrates the behavior of the control input corresponding to the interpolant operators.

	\begin{table}[ht!]
		\centering
		\caption{ Number of controllers (NC) corresponding to the interpolant operators.}
		\begin{tabular}{ c c c  }
			\toprule
			NC & Interpolant operators      \\
			\midrule
			$6$ &   Fourier modes               \\
		
			$5$&  nodal values      \\
			
			$5$ &   finite volumes        \\
		   \bottomrule
			
		\end{tabular}
	\label{table:3}
\end{table}
       

\section*{Acknowledgments}
Sudeep Kundu gratefully acknowledges the support of the Science \& Engineering Research Board (SERB), Government of India, under the Start-up Research Grant, Project No. SRG/2022/000360.
\section*{Declarations}
\textbf{Data Availability.} 

The codes  are available from authors on reasonable request. 

\textbf{CONFlLICT OF INTEREST.} 

The authors declare no conflict of interest.

\section *{Appendix }
\label{A}
 The proof of Theorem \ref{th2.1} is divided into four parts for stabilization in \(L^{2}, \ H^{1}\), and \(H^{2}\)-norms.
 
\textbf{Part 1.}
	Set \( \chi =y \) in  \eqref{1.9} to obtain
	\begin{align*}
		\frac{1}{2}	\frac{d}{dt}\norm{y}^{2}+ \nu \norm{y_{x}}^{2}-\gamma \norm{y}^{2}+\delta \norm{y}_{L^{4}}^{4}=-\mu (I_{h}(y),y)=-\mu\big(I_{h}(y)-y,y\big)-\mu\norm{y}^{2}.
	\end{align*}
	Applying the Cauchy-Schwarz and Young's inequality to the right-hand side of the above equation yields
	\begin{align*}
		\frac{1}{2}	\frac{d}{dt}\norm{y}^{2}+ \nu \norm{y_{x}}^{2}-\gamma \norm{y}^{2}+\delta \norm{y}_{L^{4}}^{4}\leq\frac{\mu}{2}\norm{I_{h}(y)-y}^{2}-\frac{\mu}{2}\norm{y}^{2}.
	\end{align*}
	Using \eqref{1.8}, it follows that
	\begin{align*}
		\frac{1}{2}	\frac{d}{dt}\norm{y}^{2}+ \nu \norm{y_{x}}^{2}-\gamma \norm{y}^{2}+\delta \norm{y}_{L^{4}}^{4}\leq\frac{\mu c_{p}^{2}h^{2}}{2}\norm{y}_{1}^{2}-\frac{\mu}{2}\norm{y}^{2}.
	\end{align*}
Rewrite the above equation as 
	\begin{align*}
	\frac{1}{2}	\frac{d}{dt}\norm{y}^{2}+ \nu \norm{y_{x}}^{2}-\gamma \norm{y}^{2}+\delta \norm{y}_{L^{4}}^{4}\leq\frac{\mu c_{p}^{2}h^{2}}{2}\norm{y_{x}}^{2}-\frac{1}{2}(\mu-\mu c_{p}^{2}h^{2})\norm{y}^{2}.
\end{align*}
Therefore, multiplying by \( 2e^{2\alpha t} \) in the above inequality and applying \eqref{2.1} and \eqref{2.2}, we observe that
\begin{align*}
 	\frac{d}{dt}(\norm{y}^{2}e^{2\alpha t})+ e^{2\alpha t}\left(\beta  \norm{y_{x}}^{2}+  2\delta \norm{y}_{L^{4}}^{4}\right)+\beta_{1}e^{2\alpha t}\norm{y}^{2}\leq 0,
 \end{align*}
where \( 0\leq \beta=2\nu- \mu c_{p}^{2}h^{2}, \) and \( 0\leq\beta_{1}=(\mu- 2\gamma-2\nu-2\alpha) \).

Integrating with respect to the time from \(0 \) to \(t \) and  multiplying by \(e^{-2\alpha t} \) in the resulting inequality, the proof of the part \(1\) is completed.

\textbf{Part 2.}
	Using the \(L^{2} \)-inner product between \eqref{1.4} and \(-y_{xx} \), we arrive at
	\begin{align*}
		\label{2.4}
			\frac{1}{2}	\frac{d}{dt}\norm{y_{x}}^{2}+ \nu \norm{y_{xx}}^{2}-\gamma \norm{y_{x}}^{2}+ 3\delta \int_{0}^{1} y^{2}y_{x}^{2}dx =\mu (I_{h}(y)-y,y_{xx})+\mu(y, y_{xx}).
		\end{align*}
Using the Young's inequality with \eqref{1.8} and multiplying  by \( 2e^{2\alpha t} \). Therefore, using \eqref{2.2} and \eqref{2.1} , we finally arrive at the following inequality
\begin{align*}
	\frac{d}{dt}(\norm{y_{x}}^{2}e^{2\alpha t})+ \nu e^{2\alpha t} \norm{y_{xx}}^{2}+6\delta e^{2\alpha t}\norm{yy_{x}}^{2} \leq 2\mu e^{2\alpha t}\norm{y}^{2}+ \mu e^{2\alpha t}\norm{y_{x}}^{2}.
\end{align*}
Integrating the above inequality over \([0, t]\) and multiplying by \(e^{-2\alpha t}\), we get from Part \(1\)
\begin{align*}
	\norm{y_{x}}^{2}+\nu e^{-2\alpha t}\int_{0}^{t}e^{2\alpha s}\norm{y_{xx}(s)}^{2}ds+6\delta e^{-2\alpha t}\int_{0}^{t}e^{2\alpha s}\norm{yy_{x}(s)}^{2}ds\leq C\norm{y_{0}}_{1}^{2}e^{-2\alpha t}.
\end{align*}
This completes the part \(2\) of the proof.

\textbf{Part 3.}
		Substituting $\chi = y_{t}$ into \eqref{1.9}, we obtain
		\begin{align*}
				\norm{y_{t}}^{2}+\frac{\nu}{2}\frac{d}{dt}\norm{y_{x}}^{2}-\frac{\gamma}{2}\frac{d}{dt}\norm{y}^{2}+\frac{\delta}{4}\frac{d}{dt}\norm{y}_{L^{4}}^{4}=-\mu (I_{h}(y)-y,y_{t})-\mu (y,y_{t}).
	\end{align*}
		Following the same step of Part \(1-2\) with \eqref{2.2}, it follows that
		\begin{align*}
				e^{2\alpha t}\left(\nu\norm{y_{x}}^{2}+(\gamma+2\nu)\norm{y}^{2}+\frac{\delta}{2}\norm{y}_{L^{4}}^{4}\right)
				+\int_{0}^{t} e^{2\alpha s}\norm{y_{t}(s)}^{2}ds
				\leq C (\norm{y_{0}}_{1}^{2}+ \norm{y_{0}}_{L^{4}}^{4}).
			\end{align*}
		The proof of the part \(3\) is completed.
		
\textbf{Part 4.}
Differentiating \eqref{1.9} with respect to \( t \) and setting \( \chi = y_{t} \), we obtain
\begin{align*}
	\frac{1}{2} \frac{d}{dt} \norm{y_{t}}^{2} + \nu \norm{y_{xt}}^{2} + 3\delta \int_{0}^{1} y^{2} y_{t}^{2} \,dx = -\mu \left(\frac{\partial}{\partial t} I_{h}(y), y_{t} \right) + \gamma \norm{y_{t}}^{2}.
\end{align*}
Applying Young's inequality to the first term on the right-hand side of the above equation, we arrive at
 \begin{align*}
	 	\frac{1}{2} \frac{d}{dt}\norm{y_{t}}^{2}+(\nu-\frac{c_{p}^{2}h^{2}\mu}{2})\norm{y_{xt}}^{2}+3\delta \int_{0}^{1}y^{2}y_{t}^{2} dx \leq (\gamma-\mu +\frac{\mu}{2}+\frac{c_{p}^{2}h^{2}\mu}{2}) \norm{y_{t}}^{2}.
 \end{align*}
Multiplying the above inequality by \( 2e^{2\alpha t} \) and integrating over the interval \( [0, t] \), we obtain a weighted energy estimate. Subsequently, multiplying the resulting expression by \( e^{-2\alpha t} \) allows us to express the inequality in a more convenient form using \eqref{2.2}
\begin{align*}
	\nonumber \norm{y_{t}}^{2}+\beta e^{-2\alpha t}\int_{0}^{t} e^{2\alpha s} \norm{y_{xt}(s)}^{2}ds+6\delta e^{-2\alpha t} \int_{0}^{t}e^{2\alpha s}(y(s)^{2},  y_{t}(s)^{2}) ds&\leq 2\alpha e^{-2\alpha t}\int_{0}^{t}e^{2\alpha s} \norm{y_{t}(s)}^{2}ds\\& \ \ +  e^{-2\alpha t}\norm{y_{t}(0)}^{2},
\end{align*}
where $ 0\leq\beta = (2\nu-c_{p}^{2}h^{2}\mu).$

From \eqref{1.4} to get
 \begin{align*}
	 		\norm{y_{t}}^{2}\leq 4\nu ^{2}\norm{y_{xx}}^{2}+4(\gamma^{2}+\mu^{2}) \norm{y}^{2}+4\mu c_{p}^{2}h^{2}\norm{y}_{1}^{2}+4\delta^{2}\norm{y}_{L^{6}}^{6}.
 \end{align*}
 Substituting $ t=0 $ to above inequality and applying Part $3$ along with \eqref{2.2}-\eqref{2.1}, we arrive at
\begin{align*}
	\norm{y_{t}}^{2}+ \beta e^{-2\alpha t}\int_{0}^{t}\norm{y_{xt}(s)}^{2}ds \leq Ce^{-2\alpha t} (\norm{y_{0}}_{2}^{2}+\norm{y_{0}}^{2}+\norm{y_{0}}_{L^{6}}^{6}).
\end{align*}
The proof of the part \(4\) is completed.
Again using \eqref{1.4} with the Sobolev embedding and applying  Part $1-4$, we arrive at
\begin{align*}
\nu \norm{y_{xx}}^{2}\leq Ce^{-2\alpha t} (\norm{y_{0}}_{2}^{2}+\norm{y_{0}}_{L^{4}}^{4}+\norm{y_{0}}_{L^{6}}^{6}).
\end{align*}
In a similar fashion, we complete the proof of Theorem~\ref{th2.1} by differentiating~\eqref{1.4} with respect to time, and then taking \(L^{2}-\)inner product with \(-y_{xxt}\), and \(y_{tt}\). In this step, we require the regularity assumption \(y_{0}\in H^{3}(0,1)\).

\end{document}